\documentclass[runningheads]{llncs}
\usepackage[T1]{fontenc}
\usepackage{graphicx}
\usepackage{amssymb}
\usepackage{amsfonts}
\usepackage{mathtools}
\usepackage{comment}
\usepackage{dblfloatfix}
\usepackage{graphics}
\usepackage{pstricks,pst-node}
\usepackage{tikz}
\usepackage{graphicx}
\usepackage[absolute]{textpos}
\usepackage{colortbl}
\usepackage{array}
\usepackage{wrapfig}
\usepackage{dblfloatfix}
\usepackage{caption}
\usepackage{pdfpages}
\usepackage{hyperref}

\newtheorem{prop}{Proposition}
\DeclareMathOperator{\Span}{span}

\begin{document}

\title{Random Schreier graphs as expanders\thanks{A preliminary result of this article has been published in the proceedings of the 30th International Computing and Combinatorics Conference (COCOON 2024). The present paper is a preprinted version of \href{https://doi.org/10.1007/s10878-026-01445-0}{doi.10.1007/s10878-026-01445-0}}}


\author{Geoffroy Caillat-Grenier}

\institute{LIRMM, Univ Montpellier, CNRS, Montpellier, France
\email{geoffroy.caillat-grenier@lirmm.fr}}

\authorrunning{Geoffroy Caillat-Grenier}      %

\maketitle

\begin{abstract}
Expander graphs, due to their mixing properties, are useful in many algorithms and combinatorial constructions. One can produce an expander graph with high probability by taking a random graph (e.g., the union of $d$ random bijections for a bipartite graph of degree $d$). 
This construction is much simpler than all known explicit constructions of expanders and gives graphs with good mixing properties (small second largest eigenvalue) with high probability. However, from the practical viewpoint, it uses too many random bits. The natural idea is to restrict the class of the bijections that we use. For example, if both sides are vector spaces $\mathbb{F}_q^k$ over a finite field $\mathbb{F}_q$, we may consider only \emph{linear} bijections, making the number of random bits polynomial in $k$ (and not $q^k$).

In this paper we provide some experimental data that shows that this approach conserves the mixing properties (the second eigenvalue) for several types of graphs (undirected regular and biregular bipartite graphs). We also prove some upper bounds for the second eigenvalue (though they are quite weak compared with the experimental results) and a lower bound for a special class of graphs.
Finally, we discuss the possibility to decrease the number of random bits further by using Toeplitz matrices; our experiments show that this change makes the mixing properties only marginally worse while the number of random bits decreases significantly.

\end{abstract}

\section{Introduction}

Expander regular graphs have been extensively studied for several decades and applications can be found in a wide variety of fields. It is remarkable that the expansion properties, which are of combinatorial nature, are connected to some spectral properties of the adjacency matrix of the graph, namely the gap between the two largest magnitude eigenvalues (the spectral gap). A large spectral gap means ``good" expansion properties. We refer the reader to \cite{Expander} which is a very complete survey on the subject. We will deal here only with \emph{undirected graphs}.

Intuitively, a graph with a large expansion ratio is ``well'' connected, and this is a desirable property in many applications, in particular in coding theory (see e.g. \cite{sipserSpielman},  \cite{zemor} or more generally \cite{Expander}, section 12). Moreover, several models used to obtain such graphs are deeply connected to mathematical structure that can be studied through the lens of graph theory, (for instance, Cayley or Schreier graphs). The expansion properties are often hard to manipulate or compute in practice, and this is where the connection with the spectrum of the adjacency matrix becomes the preferred way of dealing with expanders.

Expansion properties of graphs are interesting also because they consist of a subtle trade-off between this connectivity and the edge density of the graph. A graph that has very high edge density will be trivially an expander (e.g. the complete graph) but this is not very useful or interesting for applications or as a mathematical object. The constraint for an expander graph to be sparse, namely, to have relatively few edges (such as a \emph{logarithmic times linear} number of edges in the number of vertices i.e., with $O(\log n)$ edges incident to each of $n$ vertices) makes them much trickier to obtain. Several studies were conducted in the past decades (from the 70's) to generate such objects. It appears that \emph{random} graphs can be a particularly relevant choice to obtain expanders.

We will focus here on regular graphs. In the spectral side of things, we have some very important results on the spectral gap. The Alon-Boppana bound refers to the asymptotic lower bound for the second largest eigenvalue of regular graphs of degree $d$: $|\lambda_2| \ge 2\sqrt{d - 1} - \epsilon$ (\cite{alon-boppana}). A famous result due to Friedman (\cite{Friedman}) states that almost all such regular graphs have second largest eigenvalue asymptotically smaller than $2\sqrt{d - 1} + \epsilon$. This upper bound matches (in a probabilistic way) the lower bound which makes it optimal in some sense. Thus, most graphs have second largest eigenvalue as small as one can hope for, which explains the relevance of the random choice. 

Bipartite graphs are also very interesting and useful (in particular in coding theory). The results mentioned above do have their counterparts for such objects: we have a matching upper (\cite{bipartitBorneSup}) and lower (\cite{BipartiBornInf}) bound on the second largest eigenvalue of bipartite biregular graphs, and its nature is similar to its analog for regular non-bipartite graphs. For this reason, it is natural to use randomness to obtain bipartite spectral expanders.

The aforementioned results suggest that one can obtain and use in practical applications expander graphs by choosing them at random. However, producing such a graph, keeping it in memory and manipulations with it may require much resources. Explicit constructions of expanders (see e.g. \cite{explicit0}, \cite{explicit1}, \cite{explicit2}, \cite{explicit3}, \cite{explicit4}, \cite{explicit5}) may seem more attractive since we do not need to keep such a graph in memory: neighbors of each vertex can be computed from scratch when needed. However, this approach has its own disadvantages. Known explicit constructions of graph of expanders may lack optimality of parameters. They often imply constraints on the degree and the implementation can be difficult in practice. Besides, “polynomial time” algorithms involved in these constructions may be good in theory but too resource-demanding in practice.

The advantage of the random graphs is the simplicity: we generate a graph using a simple randomized algorithm; if it is feasible, one can compute the second largest eigenvalue to check whether the graph is suitable for applications. In some cases the graph might be too large for this calculation. Because of this, it is useful to have theoretical guarantees on the spectrum of random graphs. The Friedman bound mentioned above gives such result. However, it applies to a class of graph that can be arguably hard to sample, since the need for randomness is very high. It is then natural to try intermediate approaches (one could speak of ``pseudo-random'' models of graphs). Such graphs can be obtained using the \emph{Schreier graphs} construction with the appropriate group. These graphs are defined by a random subset of fixed size of the given group. In this work, we will study both experimentally and theoretically expander graphs, in particular the latter model.

\begin{figure}[!ht]
\centering
\includegraphics[trim={10 23 5 21},clip,scale = 0.7]{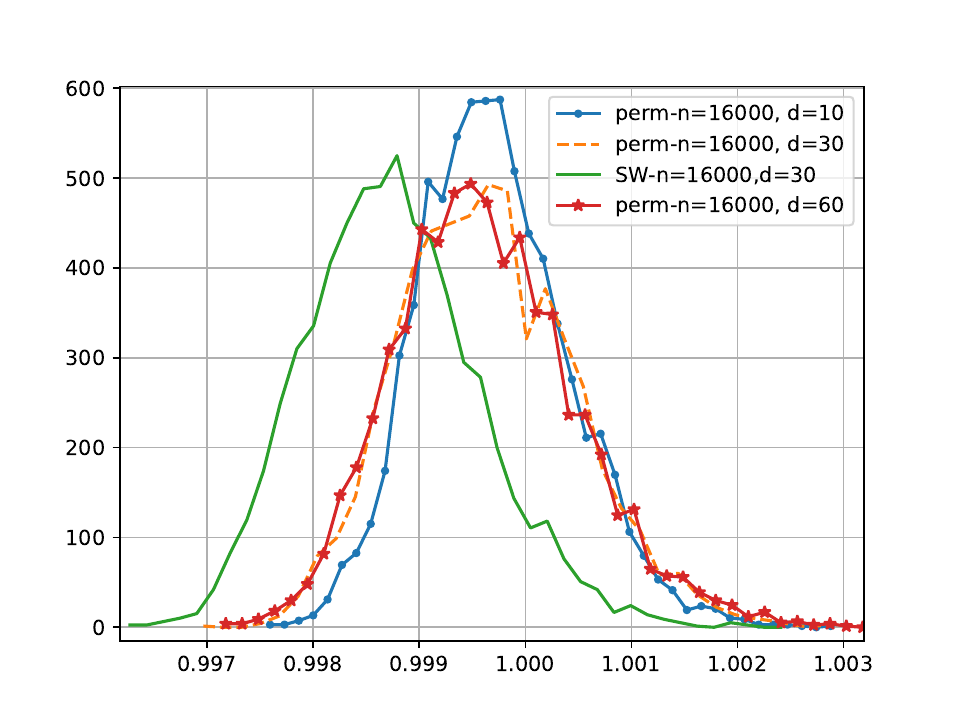}
\caption{All the eigenvalues (horizontal axis) are divided by $2\sqrt{d-1}$. The green plain curve is obtained from the uniform model \cite{unifModel}. The three others are from the permutation model for graphs of size $n = 16000$ and of degree $d=30$. Both distributions were measured by generating $5000$ graphs.}
\label{f:SW-perm}

\end{figure}

\paragraph{\textbf{Paper content.}} We begin our paper by defining all the relevant models for random graphs that will be studied (uniform model, permutation model, Schreier and Cayley graphs and Toeplitz matrices graphs). At the same time, we provide some experimental data to compare the classical models (permutation model and uniform model, Figures \ref{f:SW} and \ref{f:SW-perm}). In the next section, we will dive deeper into the experiments by computing the second largest eigenvalues of Schreier graphs of the general linear group (over a finite field) and compare it to the permutation model (truly random graphs) which will be our reference for random constructions. This is done for both non-bipartite regular graphs (Subsection \ref{RGsec}) and bipartite biregular graphs (Subsection \ref{BGsec}). Afterwards, we will present the numerical experiments for graphs from Toeplitz matrices (Schreier graphs from the general linear group restricted to the matrices that are Toeplitz, see Subsection \ref{TPsec}) and a few other spectrum measurements (on the whole spectrum) such as eigenvalues or spacings distributions (Subsection \ref{Othersec}).

Next section is dedicated to spectral theoretical bounds. We provide first a rather weak bound (Proposition \ref{WeakBound}) that applies to a very large class of graphs (it can be easily deduced from \cite{alon-roichman}). It relies on the trace method. Extending this analysis, we state and prove tighter bounds for Schreier graphs of the general linear group over the finite field of size 2 for non-bipartite (see Theorem \ref{BoundReg}) and bipartite (Theorem \ref{BoundBP}) regular graphs. Those bounds are too complicated for an asymptotic analysis, but they can give useful guarantees for specific parameters (see Tables \ref{fig:bounds} and \ref{tabBP}). We follow by proving a spectral property about our construction of biregular bipartite graphs (Proposition \ref{THbireg}). Finally, we provide some experimental data for Schreier graphs from Abelian groups (Figures \ref{fig:ab-sch} and \ref{fig:abelian}) and prove that logarithmic degree can lead to arbitrarily small spectral gap (Proposition \ref{ABstatement}; this is an alternative proof of a known result for Cayley graphs, see e.g. \cite{alon-roichman}).

\section{Models for random graphs}

We start by describing some of the main models of random graphs. The interested reader can refer to the excellent survey by Wormald \cite{wormaldModels}. Since we deal only with regular graphs, we will not discuss models such as Erdős–Rényi \cite{erdos-renyi}

\subsection{Uniform model}

Let us first present the uniform model (see section 2 in \cite{wormaldModels}). This model will not be studied deeply, but we expose some of the main properties for context. This is a model for random \emph{simple} graphs, that is, the adjacency matrix is binary and there are only zeros on the diagonal (in other words, there is at most one edge per pair of vertices and no self loops). The idea is simply to consider ${\cal G}_{n,d}$ the set of simple $d$-regular graphs (counting also the isomorphisms) with $n$ vertices (with even $dn$) and to take uniformly at random a graph of this class. 

One could observe that such a model is highly impractical. Indeed, sampling a random graph from the set ${\cal G}_{n,d}$ would require large amounts of resources (e.g. by defining a bijection between a set of integers and ${\cal G}_{n,d}$ and then producing the $i$-th graph). In \cite{unifModel}, a much more practical algorithm is presented. It was shown that the probability distribution of its outputs converges with $n$ to that of the uniform model for $d$ up to $O(n^{1/28})$. Later, in \cite{unifModelImproved}, it was proven that this distribution still converges to the desired one with $d \le O(n^{1/3})$.

\begin{figure}[ht]
\begin{center}
\includegraphics[trim={32 23 45 21},clip,scale = 0.7]{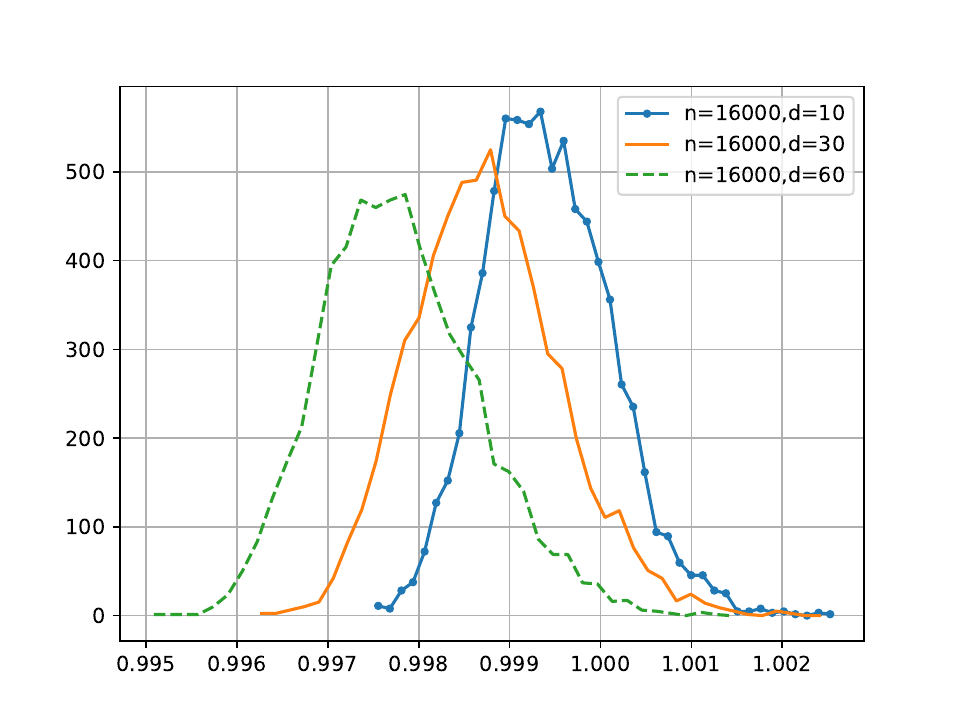}
\caption{All the eigenvalues (horizontal axis) are divided by $2\sqrt{d-1}$. The curves show the empirical distributions for the algorithm from \cite{unifModel} for graphs of size $n = 16000$ and of degree $d=10$, $d=30$ and $d=60$. These distributions were measured by generating $5000$ graphs for each.}
\label{f:SW}
\end{center}
\end{figure}

The algorithm consists in producing one random perfect matching of $nd$ vertices (this is possible since $nd$ is even) and then merging the vertices by groups of $d$. One needs to check each time an edge is selected, that it does not create loops or multiple edges with the previously chosen edges. This is rather simple to implement and the eigenvalues of the produced graphs can be measured. The empirical probability distribution of their second largest eigenvalue is displayed in Figure \ref{f:SW} for a few values of the degree.

One can observe that most of graphs given by this algorithm are \emph{Ramanujan} (i.e. their second largest eigenvalue is smaller than $2\sqrt{d-1}$). Moreover, the measured distribution seem to be a bit sensitive to the degree of the graph (at least with this algorithm). This is not too surprising since graphs of larger degree tend to have larger spectral gap (think of the complete graph which maximizes the spectral gap). This gives strong guarantees on the expansion ratio and connectivity properties of the graph. We can keep this figure in mind for what comes next in order to have insight of what can be expected from an expander graph in terms of eigenvalue. 

It is important to notice the drawbacks of this model. The algorithm executed to produce these graphs uses a lot of random bits (for each edge sampled we need $\log n$ random bits and there are $\frac{nd}{2}$ edges; the expected number of trials is $O(1)$ in most of cases, hence the space complexity is $O(nd\log n)$). Besides, it is often desirable to be able to compute all the neighbors of a given vertex (e.g. to execute a random walk in the graph) without producing and storing the whole graph. With the uniform model, there is no easy way of achieving this property. 

We would like to simplify the model in order to find a way to use less randomness in the production of such graphs. In what follows, the impact of this modification on the eigenvalues will be studied.

\subsection{Permutation model}

Our study is followed by exposing a much simpler model for random graphs that can be much more easily studied and modified. As we will see, the following model uses roughly the same amount of random bits as the uniform model. However, the structure of the graph is conceptually simpler to store and manipulate. These advantages come with the ``sacrifice'' of two properties that graphs from the uniform model have: \emph{simplicity} (i.e. there are no multiple edges) and \emph{no self loops}. Let us proceed with the definition.

\begin{definition}[Permutation model]
Select randomly and independently $d$ elements from the symmetric group $S_n$ (the group of permutations of the set $V=\{1, \dots, n\}$). This forms the multiset 
\[
S =  \{\pi_1, \dots, \pi_d\}.
\]
A $2d$-regular graph on $n$ vertices in the permutation model is then constructed as follows:
connect every vertex i with the vertex $\pi_j (i)$ for every $j \in [\![ 1, d ]\!]$. We obtain the graph 
\[
\begin{array}{lll}
&G = (V, E) \text{ with}\\
&V=\{1, \dots, n\} \text{ and}\\
&E = \{ (i,j): i = \pi(j), \pi \in S \}.
\end{array}
\]
\end{definition}

Since the graph obtained is undirected, every vertex $i$ will be connected to
$\pi_j(i)$ and $\pi^{-1}_j(i)$ for all $j \le d$ (indeed, the vertex $\pi^{-1}_j(i)$ is connected to $\pi_j(\pi^{-1}_j(i)) = i$). It also can have multiple edges and self loops. 

In this model, it is known that most of graphs have a nearly optimal spectral gap (in \cite{Friedman} it is proven that almost all graphs have second largest magnitude eigenvalue smaller than $2\sqrt{2d - 1} + \epsilon$). Thus, using this construction should allow one to get an expander with high probability.

However, as simple as it may seem, using the permutation model for sampling graphs is resource costing, particularly when it comes to store them. Observe that a randomly chosen permutation $\pi_i$ requires $\log(n!)$ bits to keep it in the memory of a computer. Hence, a description of a random graph sampled in the permutation model occupies $\Theta(d n \log n)$ bits of information (substantially, we have to store for each vertex the list of its neighbors, and there is no way to compress such a naive representation). On the other hand, it is often useful to be able to get all the neighbors of a vertex without generating the whole graph, which makes the use of randomness more limited. 

To arrive at a more educated guess on what kind of graphs such model can yield, we can, as in the previous subsection, give a histogram of the second largest eigenvalue of graphs from the permutation model. The second largest eigenvalue distribution for graphs from the uniform model is also displayed for comparison. Unlike the uniform model, for the permutation model, the degree of the graph does not seem to significantly affect the second largest eigenvalue distribution. This can be explained by the fact that the graphs obtained are not simple.

\subsection{Schreier and Cayley graphs}

In what follows we deal with a more restrictive model of random graphs (which seems to be more convenient from the practical viewpoint). Let us proceed with the definition of a \emph{Schreier graph}. 

\begin{definition}[Random Schreier graphs]
Let $H$ be a group acting transitively on a set $X$. We denote $h.x \in X$ the product of the action of $h \in H$ on $x \in X$. Let $S$ be a random multiset of elements of $H$ of $d$ elements.
Then the Schreier (undirected) graph 
$$
G = \textit{Sch}(H \circlearrowleft X, S)
$$
is a $2d$-regular graph of size $|X|$ whose edges are $(x, s.x)$ for $x \in X$ and $s \in S$.
\end{definition}

\begin{remark}
    Graphs from the permutation model can be seen as Schreier graphs of the symmetric group $S_n$ acting on $\{1, \dots , n \}$.  
\end{remark}

The graphs in our construction are Schreier graphs $Sch(H \circlearrowleft X, S)$ where $X=(\mathbb{F}_q^k)^*$ (with prime $q$; it is the vector space of dimension $k$ over $\mathbb{F}_q$ without the zero vector), $H = GL_k(\mathbb{F}_q)$ acting on $X$ by the usual matrix-vector multiplication, and $S$ is a randomly chosen multiset of group elements.

This approach uses much less bits than the general permutation model: for the case $q=2$, about $k^2$ bits for each random matrix are needed instead of $n\log n$ bits for a random permutation. We only need to store $d$ matrices of size $k\times k$  with $k=\log n$ (given the matrices from the set $S$ we can easily compute the neighbors of any vertex $v$ represented by a $k$-vector: we multiply this vector by matrices from $S$ and by its inverses).  Increasing the size of the field, we get additional economy. The problem is that there are no theoretical guarantees for mixing properties: results proven for the random permutation model do not have yet their counterparts for random \emph{linear} permutations (though some weaker results exist; see \cite{schreier2}, \cite{sabatini2021random} or \cite{schreier3}). Note, however, that this approach can be used in practice safely, assuming we can check the eigenvalues before using the graph. 

We can define \emph{Cayley graphs} as Schreier graphs whose group action is on the group itself: for a group $H$ and $S$ a multiset of elements of $H$ we have
\[
\textit{Cay}(H, S) = \textit{Sch}(H \circlearrowleft H, S).
\]
All results on Schreier graphs can therefore be applied to Cayley graphs. Actually, the famous paper about the second largest eigenvalue of Cayley graphs \cite{alon-roichman} will be one of our starting points for our theoretical results. It has inspired many results about both Cayley and Schreier graphs.

\paragraph{Practical implementation.} We now briefly explain how our construction can be implemented. The main task that needs to be done is obviously sampling the $d$ matrices from $\textit{GL}_k(\mathbb{F}_q)$. One way of doing so is simply picking uniformly and independently all $k^2$ coefficients of the matrix and reject the matrix if its determinant is zero (modulo $q$ since we are over fields whose size is a prime number). Computing the determinant is polynomial in $k = O(\log n)$. In order to estimate the expected value of the number of tries, we need to get an invertible matrix, we can estimate the following quantity:

$$
\frac{|\textit{GL}_k(\mathbb{F}_q)|}{|\mathcal{M}_{k,k}(\mathbb{F}_q)|} = \frac{\prod_{i=0}^{k-1}q^k-q^i}{q^{k^2}} = \prod_{i=0}^{k-1}1-q^{i-k}= \prod_{i=1}^{k}1-q^{-i}
$$
where $\mathcal{M}_{k,k}(\mathbb{F}_q)$ is the set of all $k \times k$ matrices whose coefficients are in $\mathbb{F}_q$.
One can show (See e.g. \cite{oeis})  that this product is convergent and always bigger than $\frac{1}{4}$. From this, we can conclude that on average, less than four tries are needed to get a non-singular matrix. Once our $d$ matrices are obtained, it is enough to apply them to all vectors. This can be done in $O(k^2)$ steps. This way, the whole graph can be produced in time $O(nd\log^2 n)$.

In next section, we will present numerical experiments for these graphs, showing the empirical behavior of their second largest eigenvalue. Later in this paper, we will proceed to an attempt to explain the observed result by showing an upper bound on the expected second largest eigenvalue of such graphs.

\subsection{Graphs from Toeplitz matrices}

We now modify our construction in order to reduce further the number of random bits. Instead of taking random invertible matrices, we take invertible Toeplitz matrices. Recall that a matrix $T = (t_{ij})$ is called a \emph{Toeplitz matrix} if $t_{ij}$ depends only on the difference $j-i$, i.e., if

\[
T=
\begin{pmatrix}
a_0      & a_1    & a_2    & \cdots & \cdots & \cdots & a_{k-1}\\
a_{-1}   & a_0    & a_1    & \ddots &        &        & \vdots\\
a_{-2}   & a_{-1} & a_0    & \ddots & \ddots &        & \vdots\\
\vdots   & \ddots & \ddots & \ddots & \ddots & \ddots & \vdots\\
\vdots   &        & \ddots & \ddots & a_0    & a_1    & a_2\\
\vdots   &        &        & \ddots & a_{-1} & a_0    & a_1\\
a_{-k+1} & \cdots & \cdots & \cdots & a_{-2} & a_{-1} & a_0\\
\end{pmatrix}
\]
for some $a_{-k+1},\ldots,a_0,\ldots,a_{k-1}$.

If we choose the elements of such a matrix randomly from a field of prime size $q$, then we get an invertible matrix with probability $1-1/q$ (see~\cite{Toeplitz}), so it is even easier to generate them. The main interest of the construction is to decrease the number of random bits used.
A Toeplitz matrix is defined by $2k-1$ elements of the field (instead of $k^2$ for a general matrix), so, in order to produce a graph, we use only $\Theta (d \log n)$ random bits ($\Theta (d \log^2 n)$ for the previous construction). One more technical advantage is that we can multiply a vector by a Toeplitz matrix fast since it is essentially the convolution. 

One might notice that such graphs cannot be called Schreier, since the set of non-singular Toeplitz matrices does not form a group (the inverse of a Toeplitz matrix is not necessarily Toeplitz). As we will see later, in the experiments, it is hard to anticipate the consequences of such observation on the graph. One can think of graphs from Toeplitz matrices as Schreier graphs from the general linear group but in which the distribution of the matrices is restricted to the space of non-singular Toeplitz matrices. However, this construction provides regular graphs just like the Schreier graphs of the general linear group, and theoretical differences between both constructions would need further research to be fully understood. One can see graphs from Toeplitz matrices as a particular subset of the Schreier graphs of the general linear group (as we had changed the distribution when choosing the matrices), since non-singular Toeplitz matrices belong to this group.

In the following section, we will present our experimental results for this kind of graphs. As for Schreier graphs of the general linear group, the second largest eigenvalue distribution has been measured. But for now, let us talk about bipartite graphs.

\subsection{Bipartite graphs constructions}

In many applications, bipartite biregular graphs are needed. The constructions previously presented can be adapted so that they provide such a graph. We now explain how this can be done.

Let us begin with the special case of a bipartite graph where the two parts of the graph are of the same size, and the degrees of the vertices in both parts are the same. Note first that every undirected regular graph can be converted to a bipartite regular graph (with the same degree in both parts), just by doubling the set of vertices (each vertex now is present both in the left and right part) and interpreting the edges as edges between parts. The random walk in this graph is just the random walk in the original undirected graph (combined with alternating the parts), so the mixing properties of this bipartite graph are the same as for the initial undirected regular graph.

However, in the bipartite case we do not need to use a permutation together with its inverse. Instead, we consider $d$ random bijections between parts and get a biregular graph of degree $d$. We use this construction in our experiments with bipartite graphs and in our theoretical results (the previous construction would have given the same result as for regular non-bipartite graphs). This allows us to have bipartite graphs of any (possible odd) degree. On the other hand, one can observe that this way of producing such graphs uses twice as many random matrices (or permutations) as in the regular non-bipartite case. As we will see later, this does not seem to affect significantly the experimental results. However, our theoretical bounds are for this reason stronger in the bipartite case.

Sometimes (e.g., in some constructions in coding theory) we need biregular bipartite graphs where left and right parts  are of different size. A natural way to get such a graph is to start with a bipartite graph with the same part sizes, and then merge vertices in one of the parts. Merging $s$ vertices into one, we decrease the size and increase the degree by the factor of $s$. (We assume that number of vertices is a multiple of  $s$.)

Again, we conducted numerical experiments on the second largest eigenvalue of such graphs that are presented in the following section. Moreover, we have shown some spectral upper bound on the expected value that is detailed in the last sections of this paper.

\section{Experimental results}

We now expose the experimental results obtained on the graph constructions that were presented in the last section, that are Schreier graphs of the general linear group over a finite field of prime size, and graphs from Toeplitz invertible matrices (again, over the same finite fields). We will also display experiments with bipartite biregular Schreier graphs from the general linear group.

We will start with the histograms for the second largest eigenvalues of regular non-bipartite graphs from the general linear group that will be compared with that of the permutation model. Then we will show similar picture for bipartite graphs (both regular and biregular). Afterwards, we will switch to graphs from Toeplitz matrices. At the end of this section, we will show some other measurements on the spectrum of graph that can characterize randomness of the models. 

\paragraph{\textbf{Data availability statement.}} 
The source code used to generate the experiments and figures is available at:\\
\href{https://gite.lirmm.fr/gcaillatgrenier/random-schreier-graphs-as-expanders}{gite.lirmm.fr/gcaillatgrenier/random-schreier-graphs-as-expanders}

\subsection{Second largest eigenvalue of Schreier graphs of $\textit{GL}_k(\mathbb{F}_q)$}\label{RGsec}

Consider a finite field $\mathbb{F}_q$ with $q$ ($q$ prime) elements, and the vector space $\mathbb{F}_q^k$  over $\mathbb{F}_q$. Choose randomly $d/2$ invertible matrices $T_1,\ldots,T_{d/2}$, and use the corresponding permutations of the set of non-zero tuples to get an undirected graph of degree $d$: a non-zero tuple $x\in \mathbb{F}_q^k$ has neighbors
$$
T_1 x,\ldots, T_{d/2}x, T_1^{-1}x,\ldots, T_{d/2}^{-1} x.
$$
This is a regular graph with $q^k-1$ vertices and degree $d$ (that may contain loops and parallel edges), and we compute numerically the second eigenvalue of this graph.

\begin{figure}[!ht]
\begin{center}
\includegraphics[trim={32 23 45 41},clip,scale = 0.7]{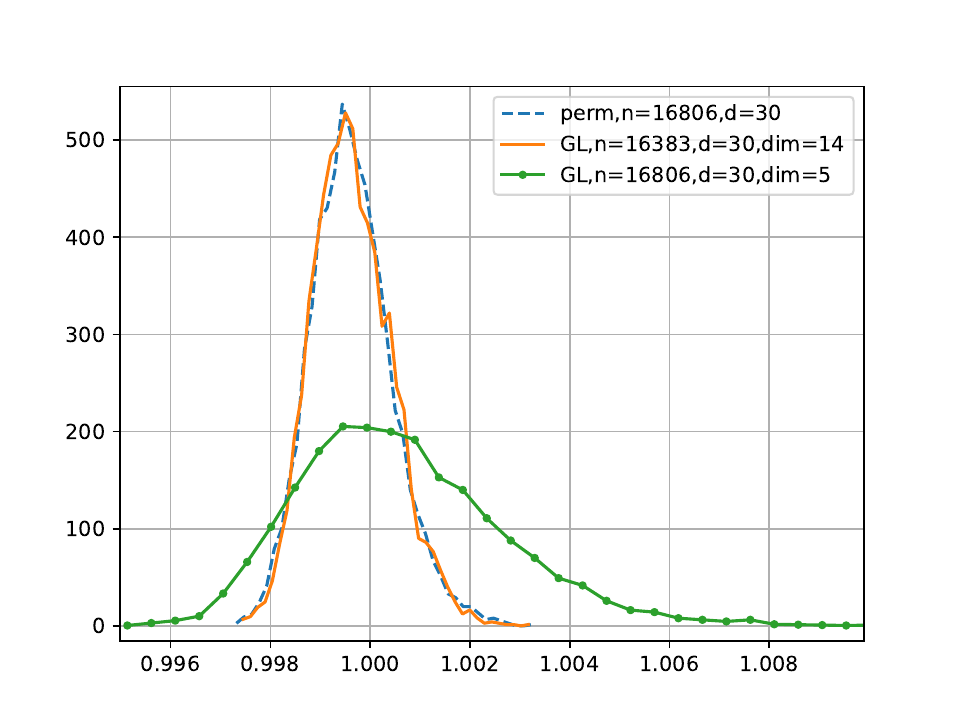}
\caption{All the eigenvalues (horizontal axis) are divided by $2\sqrt{d-1}$. The two upper curves (quite close to each other) show the empirical distributions for random permutation model (dashed line) and random \emph{linear} permutation model for graphs of degree $d=30$ ($15$ permutations) and size $16383=2^{14}-1$. The third curve (wider green one in the bottom part) shows the distribution for random linear permutations on $5$-tuples from a $7$-element field, so the graph has $7^5-1=16806$ vertices (and the same degree $30$, obtained with $15$ permutations). For each distribution $5000$ experiments were made; the eigenvalues are grouped in $40$ bins.}
\label{fig:reg}
\end{center}
\end{figure}

This procedure is repeated many times, and then the empirical distribution of the second eigenvalues is shown (Figure~\ref{fig:reg}). For comparison, we show also the distribution for the random permutation model and for a different field sizes (and approximately the same graph size).

The invertible matrices were generated by choosing a random matrix and then testing whether it is invertible. The random numbers were obtained by using random bits from a physical phenomenon (electromagnetic noise in an electronic circuit; the experiments were conducted by Alexander Shen, LIRMM, ESCAPE). As mentioned earlier, for Boolean matrices, the fraction of invertible ones is above $1/4$, so we do not need too many trials; for bigger fields, the fraction is even bigger. Computing the determinant (checking the invertibility) is easy. For computing the eigenvalues we used the \texttt{C++} library Spectra\footnote{Spectra's home page: \url{https://spectralib.org/}.} that implements the Lanczos algorithm~\cite{lanczos}.

The choice of the unit on the horizontal axis ($2\sqrt{d-1}$ for non-normalized adjacency matrices) is motivated by the asymptotic lower bounds we mentioned. Graphs that have the second eigenvalue smaller than this unit are usually called \emph{Ramanujan} graphs; we see that for our parameters both the random permutations and random linear permutations give Ramanujan graphs most of the time.

We also tested another set of parameters: $5$-tuples from a $7$-element field (so the resulting graph size is close to the previously considered one, $7^5-1=16806$ instead of $2^{14}-1=16383$, which makes the results roughly comparable). Note that this choice allows us to decrease the number of random bits used for the graph generation. We see that this is achieved for a price: the distribution is now less concentrated,  though still very close to the Ramanujan threshold. One could expect that more random bits lead to better concentration.

One can also observe (we do not have any explanation for this paradoxical observation) that the probability to get a random graph that is on the left of, say, $0.998$, becomes bigger for a bigger field. So, if we try to get one graph with small eigenvalue after many trials, the bigger field may be better. One could try to go to the extreme in this direction and consider $1\times 1$ invertible matrices, i.e., consider multiplication by non-zero elements in some field $\mathbb{F}_q$ (the graph size is then $q-1$). This is, however, a bad idea since in this case the invertible matrices commute and therefore many paths in a random walk arrive at the same place just due to commutativity. Such a graph will have poor mixing properties.

\begin{figure}[ht]
  \centering
\includegraphics[trim={30 23 30 41},clip,scale = 0.7]{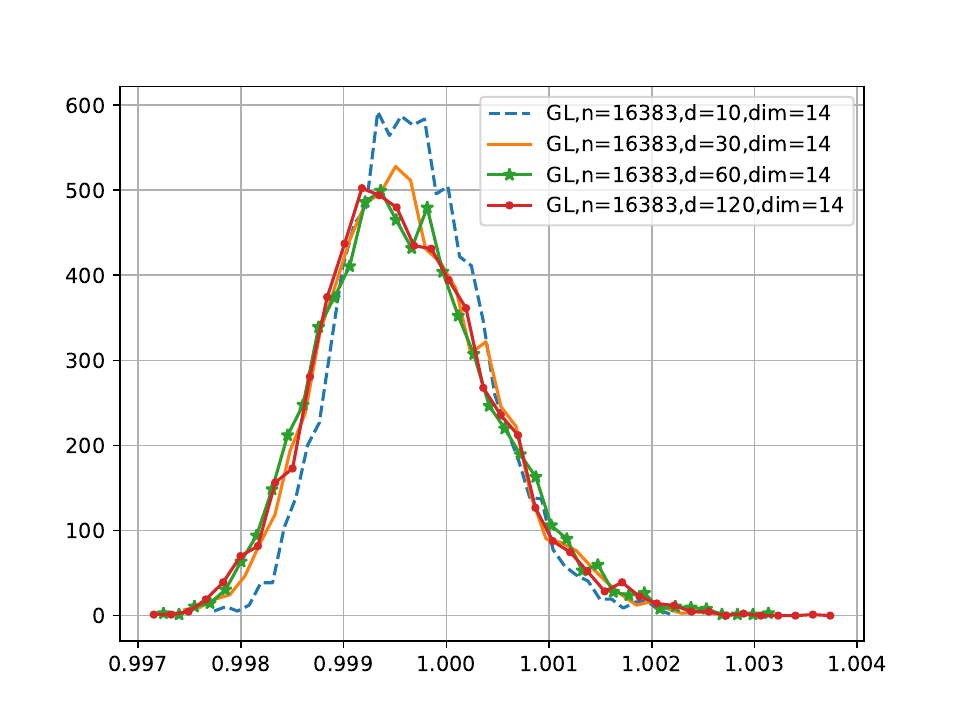}
  \caption{Second largest eigenvalue distributions for $\textit{GL}_{14}(\mathbb{F}_2)$ Schreier graphs of degrees $10, 30, 60$ and $120$.}\label{fig:reg-deg}
\end{figure}

\begin{figure}[!ht]
  \centering
    \includegraphics[trim={30 23 30 41},clip,scale = 0.7]{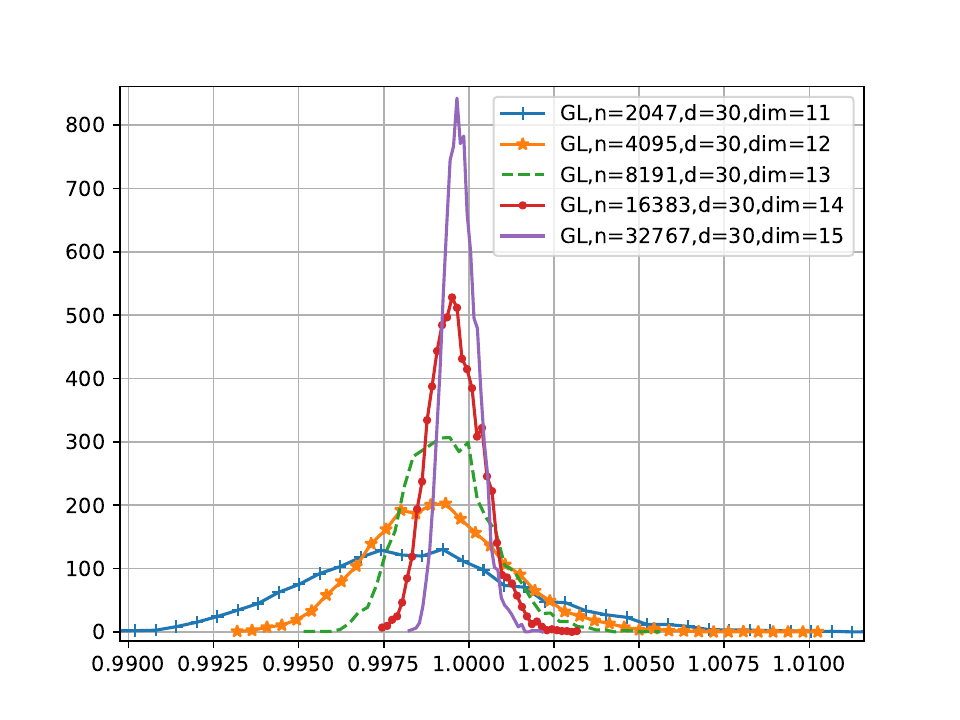}
      \caption{  Instead of the degree, we vary the dimension of the matrices and keep the same field $\mathbb{F}_2$.}\label{fig:reg-sizes}
\end{figure}

Additionally, we can observe that variations on the degree of the graphs (for a given construction) does not seem to have much impact on the distribution of the second largest eigenvalues (Figure \ref{fig:reg-deg}). Though the difference is hardly visible, we can see that smaller degree leads to a slightly greater second largest eigenvalue. This is not surprising since low degree random graphs tend to be worse expanders. However, the eigenvalues are slightly more concentrated. As expected (see Figure \ref{fig:reg-sizes}) increasing the dimension does not change the expected value but decreases the variance.

\subsection{Second largest eigenvalue of bipartite biregular graphs}\label{BGsec}

The experiments in the previous section were for undirected regular graphs, but one may also need a \emph{bipartite} graph that is \emph{biregular} (degrees of the vertices are the same in each part, but may differ in left and right parts) and has good mixing properties (for the back and forth random walk).

\begin{figure}[!ht]
\centering
\includegraphics[trim={20 23 40 41},clip,scale = 0.7]{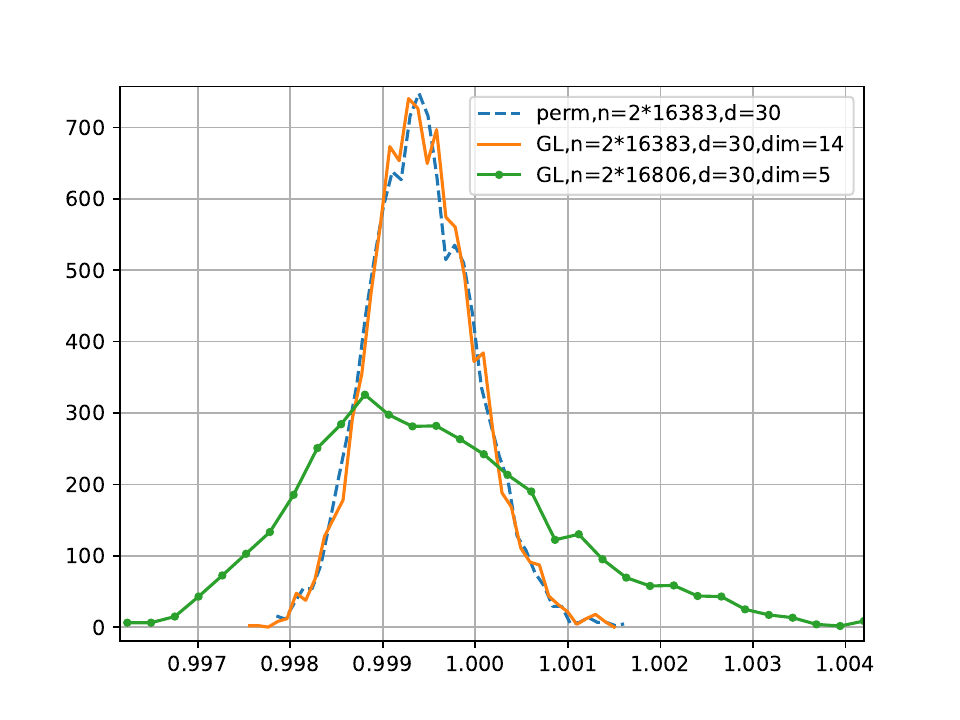}
\caption{Here, \texttt{n1} is the size of each partition. The dashed line (blue) shows the distribution of second eigenvalues for a bipartite graph constructed from $30$ random permutations of $16383$-element set (nonzero bit vectors of size $14$). The orange line (close to the first one) shows the same for $30$ \emph{linear} permutations (over $\mathbb{F}_2$). The third line (green, in the bottom part) shows the distribution for $30$ permutations of a $16806$-element set that are randomly chosen among the $\mathbb{F}_7$-linear permutations of non-zero $5$-tuples with elements in $\mathbb{F}_7$. The blue dashed curve is actually obtained using ``only'' $\approx 4600$ graphs: we ran out of random bits (10Go) for this model and these parameters.}\label{fig:BP-reg}
\end{figure}

\begin{figure}[!ht]
    \centering
\includegraphics[trim={20 23 40 41},clip,scale = 0.7]{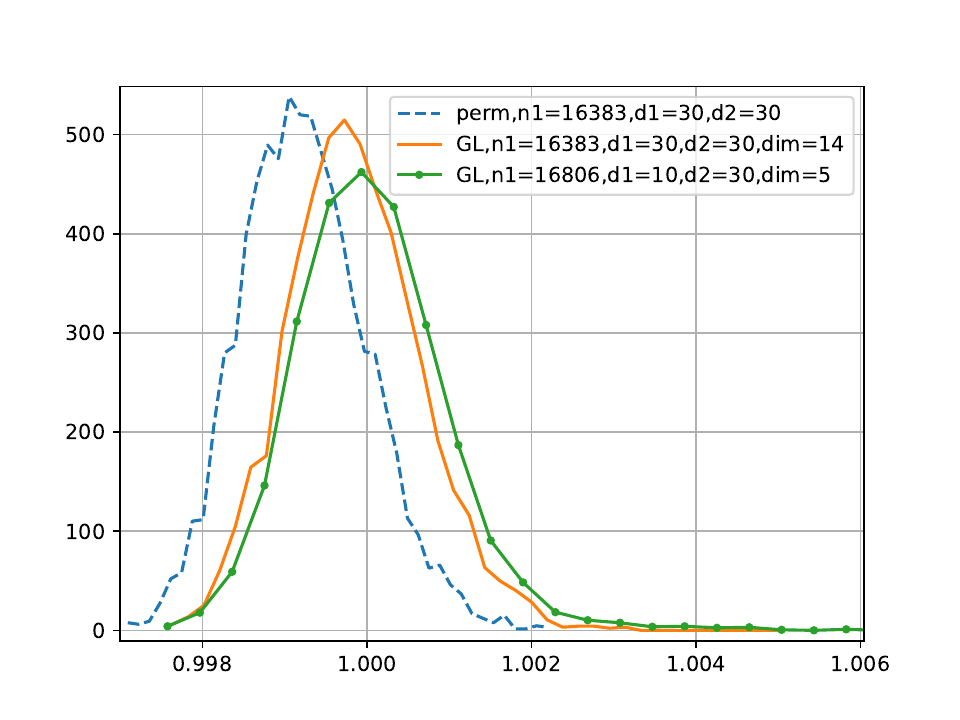}
\caption{The dashed line (blue) shows the distribution of the second eigenvalues for $5000$ experiments when we take $10$ random permutations of a set with $2^{14}-1=16383$ elements, construct a bipartite graph with degree $10$ in both sides and then merge triples of vertices on one side thus getting degree $30$ and size $16383/3=5461$. The solid (orange) line replaces random permutations by random \emph{linear} (over $\mathbb{F}_2$) permutations of non-zero vectors in $\mathbb{F}_2^{14}$. Finally, the line with dots (green) shows the same distribution if we take $10$ linear permutations of non-zero vectors in $\mathbb{F}_7^5$.}\label{fig:BP-bireg}
\end{figure}

These constructions can be performed both with random permutations and random \emph{linear} permutations. The same question arises: what can be said about mixing properties of graphs in these two settings? we recall here the relevant properties of the adjacency matrix of bipartite graphs that will help us to answer this question.

A traditional way to measure mixing property of a bipartite graph is to consider the (non-normalized) adjacency matrix for it (height and width are equal to the total number of vertices in both parts).  The eigenvalues of this matrix are grouped into pairs, starting with $\sqrt{d_Ld_R}$, where $d_L$ and $d_R$ are left and right degrees (the eigenvector coordinates are $\sqrt{d_L}$ in the left part and $\sqrt{d_R}$ in the right part) and $-\sqrt{d_Ld_R}$ (to change the sign of the eigenvalue we change signs in all eigenvector coordinates in one of the parts). The next pair, $\lambda>0$ and $-\lambda<0$, determines the mixing properties, and we call $\lambda$ the second largest eigenvalue (ignoring the negative ones).

The role of the Ramanujan threshold $2\sqrt{d-1}$ (for $d$-regular undirected graphs of degree $d$) is now played by $\sqrt{d_L - 1} + \sqrt{d_R - 1}$: it has been proven that the second largest  eigenvalue cannot be much smaller than this quantity~\cite{BipartiBornInf} and with high probability is not much larger for random graphs~\cite{bipartitBorneSup}.  We use $\sqrt{d_L-1}+\sqrt{d_R-1}$ as the unit on the horizontal axis for all our figures.

Figure~\ref{fig:BP-reg} shows the experiments with bipartite graphs with the same degree $30$ in both parts obtained using $30$ permutations of a set of non-zero elements of $\mathbb{F}_2^{14}$. The setting is different from the case of the undirected regular graph (Figure~\ref{fig:reg}) where we used $15$ permutations and their inverses, but the results turn out to be quite similar (and the distribution for the case of $\mathbb{F}_7^5$ is also similar to Figure~\ref{fig:reg}).

For the case of different left and right degrees the results of our numerical experiments are shown in Figure~\ref{fig:BP-bireg}. 
Comparing this figure with the previous one, we observe some differences: the random permutations now give slightly (but visibly) better results compared to random linear permutations. On the other hand, changing the size of the field (and matrix dimension) now has much less effect. In fact, results may depend on which triples we are merging, for the linear case,  due to the additional structure induced by the the construction (compared to general permutations). However, this phenomenon still has to be explained.

\subsection{Second largest eigenvalue of graphs from Toeplitz matrices}\label{TPsec}

We now turn back to regular graphs from Toeplitz matrices and present the experimental results we have measured on their second largest eigenvalue distribution.

\begin{figure}[!ht]
  \centering
    \includegraphics[trim={20 23 30 41},clip,scale = 0.7]{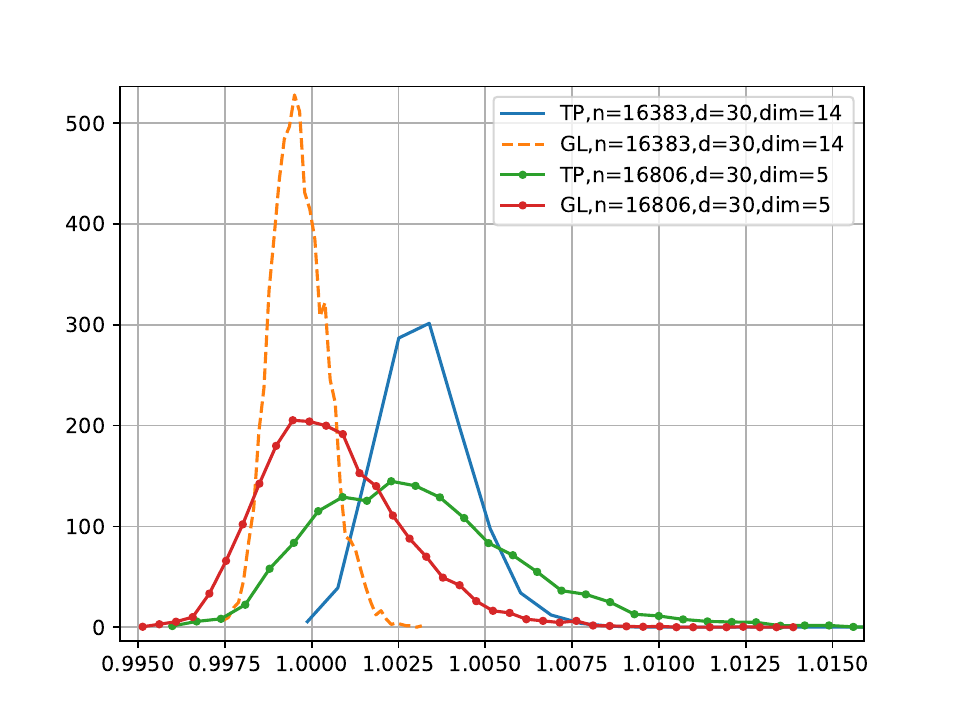}
  \caption{Distributions of second largest eigenvalues of 5000 graphs obtained from Toeplitz matrices  (\textsf{TP}) for $\mathbb{F}_2^{14}$ and $\mathbb{F}_7^5$ compared with the corresponding empirical distributions for random linear permutations in the same vector spaces (\textsf{GL}).}\label{fig:TP-GL}
  \end{figure}

  \begin{figure}[!ht]
  \centering
    \includegraphics[trim={20 23 45 41},clip,scale = 0.7]{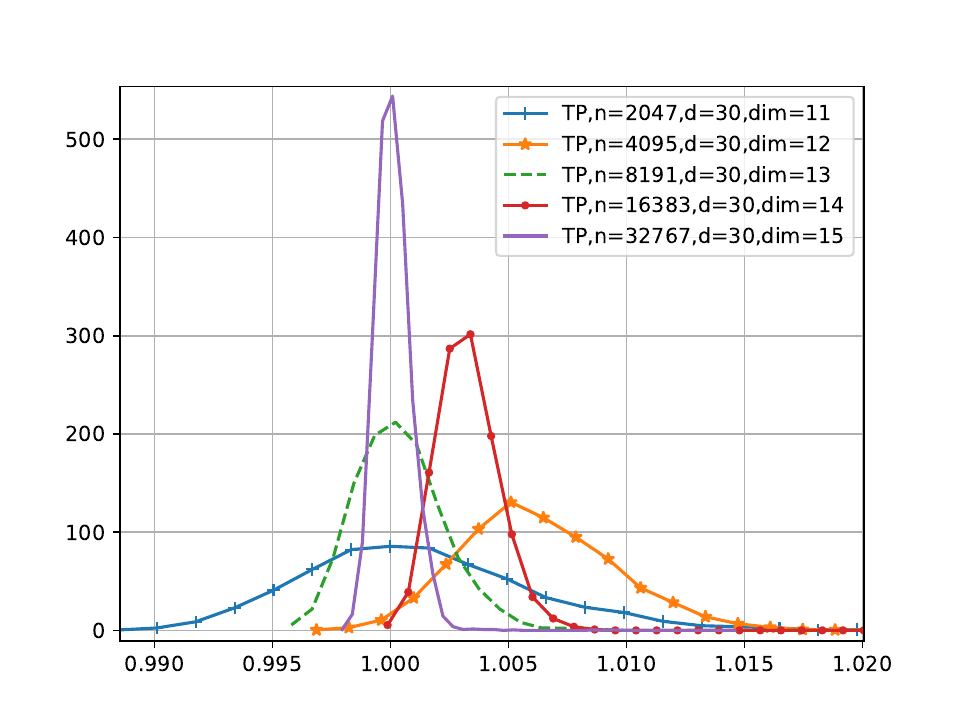}
  \caption{Here, all graphs are from Toeplitz matrices over $\mathbb{F}_2$.}\label{fig:TP-sizes}
\end{figure}

In Figure~\ref{fig:TP-GL}, we return to the case of undirected regular graphs, and compare random graphs obtained using $15$ random matrices (either random invertible matrices or random invertible Toeplitz matrices). Two upper curves (yellow and blue) use $\mathbb{F}_2^{14}$; one can see that Toeplitz matrices give bigger (and less concentrated) eigenvalues. The same can be seen from the two bottom curves (red and green); we see also that decreasing the dimension and increasing the field size increase significantly the variance (but not the mean, see below). The difference between Toeplitz and general invertible matrices is visible but small. Hence, Toeplitz matrices are suitable for producing good spectral expanders (at least with the parameters of the experiments).

In Figure \ref{fig:TP-sizes}, we observe that the parity of the dimension of Toeplitz matrices used affects the mean of the distribution (the variance behaves as expected). Odd dimensions tend to produce better spectral expanders than even dimensions (this can be also observed on preceding figure (\ref{fig:TP-GL}), red curve). We could not find any explanation for this phenomenon.

\subsection{Some other measurements on the spectrum}\label{Othersec}

There exist some other spectral properties that can be measured in random regular graphs. These properties imply, however, to get a global view on the spectrum of the graph. Therefore, in order to measure them, we need to produce significantly smaller graphs so that the computation can be done in a reasonably short amount of time. Again, we have used here the \texttt{C++} Spectra library to compute the whole spectrum of reasonably large graphs (e.g. about $3000$ vertices). 

\begin{figure}[!ht]
  \centering
    \includegraphics[trim={20 23 45 41},clip,scale = 0.6]{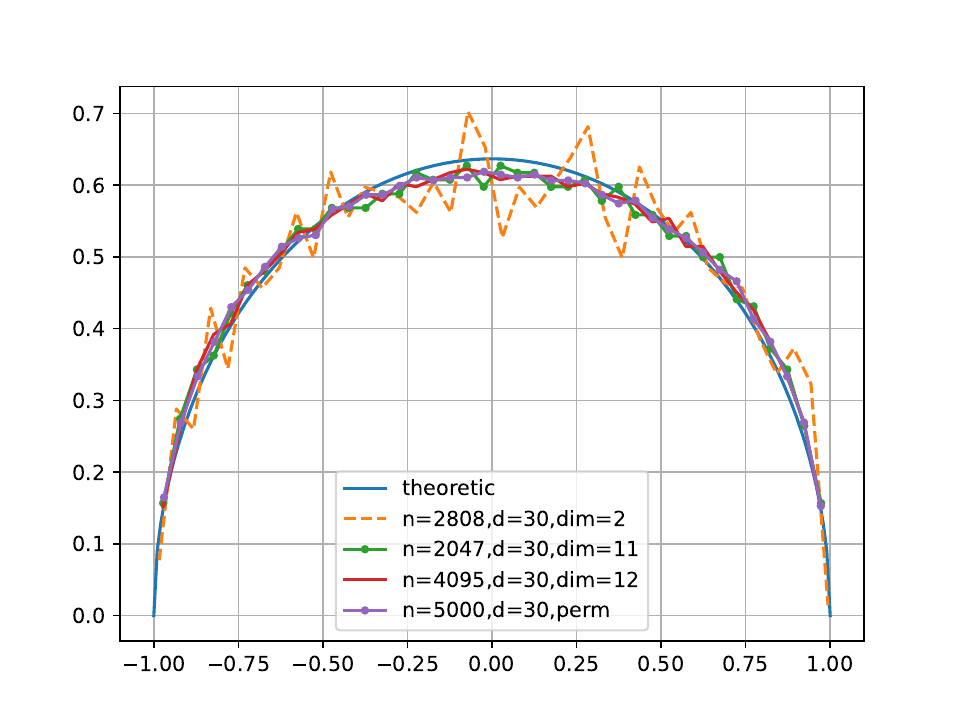}
  \caption{We display here the normalized eigenvalue histogram of four Schreier graphs. Three of them are generated with linear permutations (non-singular matrices over finite fields); the last one is from the permutation model. We can see that in large dimension, the distribution is pretty similar to that of the semicircle law. In small dimension (dimension 2), the behavior seems to be more unpredictable. We show also the theoretical distribution (\cite{spectrumDistrib}) of Wigner's semicircle law.}\label{fig:spectrum-GL}
\end{figure}

We first display in Figure \ref{fig:spectrum-GL} the distribution of the whole spectrum of the normalized adjacency matrix of a few graphs from our construction with non-singular matrices over a finite field. Such distribution is that of the Wigner's semicircle for the uniform model for random regular \emph{simple} graphs (\cite{spectrumDistrib}). We can see that such distribution seems to correspond for graphs from the permutation model, but also for graphs from our construction (Schreier graphs of the linear general group), at least in high enough dimension.

\begin{figure}[!ht]
  \centering
    \includegraphics[trim={20 23 45 41},clip,scale = 0.6]{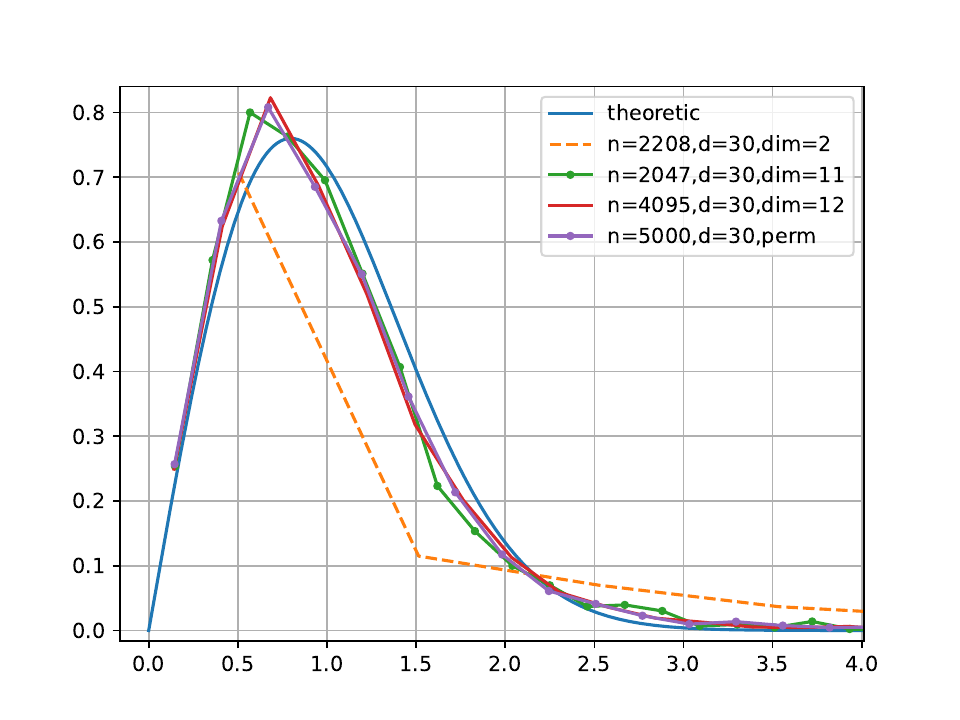}
  \caption{This figure represents the histogram of the gaps (in absolute value) between every two consecutive eigenvalues (in order of magnitude) of the spectrum of the normalized adjacency matrix of four Schreier graphs (as in the preceding figure). Again, we additionally display the curve of the distribution of spacings of the Gaussian Orthogonal Ensemble (\cite{spacings}) for comparison.}\label{fig:spacings-GL}
\end{figure}

On the other hand, it is interesting to display also the spacing distribution (i.e. the distribution of the difference between two consecutive eigenvalues). Such distribution for random regular graphs is that of the Gaussian Orthogonal Ensemble (see \cite{spacings}). One can observe that, unless the dimension  of the matrices is pretty small, Schreier graphs from the general linear group, and Schreier graphs from the symmetric group (permutation model) have spacing distribution that behaves like that of the uniform model.

It is important to notice that both figures are generated with only one graph for each curve, hence these results can be taken with a grain of salt. These experiments are here to help understanding the general behavior of random graphs from our construction and can corroborate the fact that they have similar properties to those of truly random graphs.

\section{Spectral bounds}

In this section, we will present our attempt to explain our experimental results.
As it will be noticeable below, the bounds we prove are not very explicit. In fact, we will use in most cases numerical applications of the bounds by fixing some parameters to extract more tangible facts from our messy formulas.

\paragraph{A word about the usefulness of our bounds.} This computer assisted bound can be applied to relatively large graphs. Let us mention that a ``relatively large'' graph is not necessarily a graph that can be stored in computer's memory. In practice, we may need graphs with the number of vertices exponentially larger than the available memory (e.g., we may need an expander graph with $2^{256}$ vertices, where every vertex of the graph represents a string of $32$ bytes)

One of the well-known applications of expanders is a deterministic error amplification, which helps to reduce the error probability in a randomized algorithm without increasing the number of random bits, see  \cite[Section 1.1.3]{Expander}. In applications of this type, one would need very large expander graphs, so that we cannot store the entire graph in the memory of a computer and only keep in mind the indices of a few vertices. In this case, we would need a \emph{strongly explicit} construction of an expander: given the index of a vertex, we should be able to compute efficiently the list of its neighbors without storing the whole graph.

Our constructions are strongly explicit, and  they, could be used for applications of this kind. However, we cannot compute numerically the eigenvalues of graphs whose size is exponential compared with the memory of a computer. So in this case we need theoretical bounds for the spectrum of sampled graphs.

We now proceed to prove some upper bounds on the expected second largest eigenvalue of the graphs from the construction we have presented.

\subsection{Weak bound for a large class of graphs}

It is useful to have theoretical guarantees about the eigenvalues when we cannot compute them. We give statements about the expectation of the second largest eigenvalue. The expectation is taken over all possible Schreier graphs of given size and degree of some group, i.e. among all multi-subset of appropriate size of the group considered. To do so, one can use the trace method~\cite{broder-shamir} for the case of Cayley (or arbitrary Schreier) graphs, following~\cite{alon-roichman}. 
This gives the following bound (see also \cite{sabatini2021random} and \cite{schreier2}):
\begin{prop}\label{WeakBound}
Let $G$ be a Schreier graph constructed by choosing at random $d$ elements of a group acting on an $n$-element set \textup(so it has $n$ vertices and degree $2d$\textup). Let $\mu_2$ be the second largest eigenvalue of its normalized adjacency matrix. Then
$$
\mathbf{E}(|\mu_2|) \le e\sqrt{\frac{\ln n}{d}}
$$
where $\ln$ is the natural logarithm and $e$ its base.
The same bound holds for graphs constructed with $d$ randomly chosen invertible Toeplitz matrices. The statement is also valid for bipartite graphs of $2n$ elements of degree $d$ (produced with $d$ random elements).
\end{prop}

This proposition is proven in section \ref{proofP1}. We can put this statement in the more traditional form (that can be found e.g. in \cite{alon-roichman}, \cite{schreier2}, \cite{schreier3}):

\begin{corollary}
    Let $G$ be a $2d$-regular random Schreier graph on $n$ vertices and of second largest eigenvalue $\lambda_2$. For every $\epsilon > 0$ there exists a constant $C_\epsilon$ such that
    \[
    d = C_\epsilon \ln n \implies \mathbf{E}(|\lambda_2|) \le \epsilon.
    \]
\end{corollary}

\begin{proof}
    Direct by replacing $d$ in Proposition \ref{WeakBound} with $C_\epsilon = (\frac{e}{\epsilon})^2$.
  \qed \end{proof}

\subsection{Spectral bound for regular Schreier graphs from our construction}

By focusing on Schreier graphs of $\textit{GL}_k(\mathbb{F}_2)$ and by pushing forward the analysis used in the proof of Proposition \ref{WeakBound}, it is possible to improve the previous result. The precise statement is given below (see Theorem \ref{BoundReg}). Unfortunately, the improved bound is quite complicated (and even its asymptotic behavior is not clear). However, using computer assistance, we can compute this improved bound for specific values of parameters ($k$ and $d$).

\begin{table}[!ht]
\begin{center}
   \begin{tabular}{| c || c | c | c |}
     \hline
     $(k, d)$ & Theorem \ref{BoundReg}  &  Proposition~\ref{WeakBound} & Ramanujan \\ \hline
     $(14, 15)$  & 0.7607 & 2.1863 & 0.3590\\
     $(14, 30)$  & 0.5033 & 1.5459 & 0.2560\\
     $(14, 60)$  & 0.3389 & 1.0931 & 0.1818\\
     $(20, 20)$  & 0.7758 & 2.2631 & 0.3122\\
     $(25, 50)$  & 0.4969 & 1.6002 & 0.1989\\
     $(30, 500)$ & 0.1435 & 0.5543 & 0.0632\\
     $(40, 500)$ & 0.1618 & 0.6401 & 0.0632\\
     $(50,1000)$ & 0.1217 & 0.5060 & 0.0447\\
     $(60, 60)$  & 0.7377 & 2.2631 & 0.1818\\
     $(200,200)$ & 0.7016 & 2.2631 & 0.0998\\
     \hline
   \end{tabular}
\end{center}
\caption{Upper bounds for the expected value of $\mu_2$ provided by Proposition~\ref{WeakBound} and its improved version for a few choices of parameters (the dimension $k$ of matrices of the group, and the size $d$ of the multisubset used to produce the graph, hence the degree is $2d$).  We display also the Ramanujan's threshold $\frac{2\sqrt{2d-1}}{2d}$ for comparison. The first three lines correspond to graphs that we tested in our experiments (see Figures \ref{fig:reg-deg}).}\label{fig:bounds}
\end{table}

Looking at this table, we see that the improved bound obtained with the help of Theorem \ref{BoundReg} makes sense in several cases where Proposition~\ref{WeakBound} does not give anything nontrivial (the bound is greater than the trivial bound $1$). However, even the improved bound is far from being tight (recall that our experiments, when feasible, exhibit a rather small difference between the second eigenvalue and Ramanujan's threshold). Moreover, we can observe some convergence behavior when $d=k$. This would imply that the asymptotic behavior of the improved bound is similar to that of Proposition \ref{WeakBound}.

We now can formulate the statement. To simplify the notations, here $d$ will represent the number of random matrices used. Hence the degree of the graph will be $2d$.

\begin{theorem}\label{BoundReg}
Let $G$ be a random Schreier graph of $\textit{GL}_k(\mathbb{F}_2)$ acting on $(\mathbb{F}_2^k)^*$ with the random multiset $S$, $n = 2^k - 1$ (with $k \ge 2$), and degree  $2d = 2|S|$.
 Let $\mu_2$ be the second largest eigenvalue of its normalized adjacency matrix. Consider the following recurrent relation:
\[
X_p(c,l, d) = \left\{
    \begin{array}{ll}
        1 \mbox{ if } c = 0 \mbox{ and } l = 0; \\
        0 \mbox{ if } p > l; \\
        \sum_{i = c}^{\lfloor\frac{l}{p}\rfloor} \binom{d}{i}
\bigg(\frac{2^p}{p!}\bigg)^i\frac{l!}{(l-pi)!}X_{p+1}(0,l-pi, d-i) \mbox{ otherwise.}
    \end{array} 
\right.
\]
We set $x_3(i) = X_3(1,2m-2i, d-i)$ and $x_4(i) = X_4(0,2m-2i, d-i)$.
Then, for every integer $m$, 
\begin{multline}\label{formulaTH1}
\mathbf{E}(|\mu_2|) \le 
 \bigg(\bigg(\frac{1}{2d}\bigg)^{2m}n\bigg[\sum_{i=1}^m\binom{d}{i}\frac{(2m)!}{(2m-2i)!}
\bigg[(2^i-1)(x_3(i) + x_4(i))\frac{2}{n}
\\ + \bigg((x_3(i) + x_4(i))\frac{1}{n}
 +\frac{2^i}{(i+1)!}\bigg(x_3(i)\frac{5}{n} + x_4(i) \bigg)\bigg)\bigg]
\\ + X_1(1, 2m, d)\frac{1}{n} + x_3(0)\frac{5}{n} + x_4(0)\bigg]-1\bigg)^{\frac{1}{2m}}.
\end{multline}
\end{theorem}
The bound~(\ref{formulaTH1}) looks messy, and the asymptotic analysis may be difficult. However, this formula becomes useful when we need to calculate (with the help of a computer) a bound for a specific graph of some given size.

To compute this bound for given $k$ and $d$, we choose the value of the size of the walk ($2m$) such that (\ref{formulaTH1}) implies the strongest bound for $\mathbf{E}(|\mu_2|)$. The exact value of $m$ that minimizes the bound is unknown, but we suspect it to be logarithmic as in Proposition \ref{WeakBound}. In order to compute this bound for some $k$ and $d$, we set $m=1$ and compute (\ref{formulaTH1}) right hand side several times, incrementing $m$ until it ceases decreasing (\ref{formulaTH1}) (the bound will converge to $1$ after reaching its minimum when $m$ grows). In all our computations, this happens for $m \le k$. We discovered that for  several ``reasonable'' values of $k$ and $d$, this formula implies a non-trivial bound for $\mu_2$ (see Table \ref{fig:bounds}).

\subsection{Spectral bound for bipartite graphs}

A pretty similar bound can be proved for regular bipartite graphs obtained from random non-singular binary matrices. The techniques used for the proof are the same as that of the preceding theorem.

\begin{theorem}\label{BoundBP}
Let $G$ be a random Schreier bipartite graph of $\textit{GL}_k(\mathbb{F}_2)$ ($k\ge2$) acting on $(\mathbb{F}_2^k)^*$ with the random multiset $D$, $|D| = d$, and $n=|(\mathbb{F}_2^k)^*|$. Let $\mu_2$ be the second largest eigenvalue of its normalized adjacency matrix. Consider the relation 
$$
Y_p(c,l, d) = \left\{
    \begin{array}{ll}
        1 \mbox{ if } c = 0 \mbox{ and } l = 0; \\
        0 \mbox{ if } p > l; \\
        \sum_{i = c}^{\lfloor\frac{l}{p}\rfloor} \binom{d}{i}
(p!)^{-i}\frac{l!}{(l-pi)!}Y_{p+1}(0,l-pi, d-i) \mbox{ otherwise.}
    \end{array} 
\right.
$$
We set
$y_3(i) = Y_3(1,2m-2i, d-i)$, $y_4(i) = Y_4(0,2m-2i, d-i)$. Then, for every integer $m$, 
\begin{multline}\label{formulaTH2}
\mathbf{E}(|\mu_2|) \le \bigg(\bigg(\frac{1}{d}\bigg)^{2m}n\bigg[\bigg(\frac{1}{d}\bigg)^{2m}
\\ \sum_{i = 1}^m \bigg[\binom{d}{i} \bigg(\frac{1}{2}\bigg)^i \frac{(2m!)}{(2m-2i)!}
\frac{2^i}{(i+1)!}\bigg(y_3(i)\frac{5}{n} + y_4(i) \bigg)\bigg]
\\ + Y_1(1, 2m, d)\frac{1}{n} + Y_2(1, 2m, d)\frac{2}{n}  + y_3(0)\frac{5}{n} + y_4(0)\bigg] - 1 \bigg)^{\frac{1}{2m}}.
\end{multline}

\end{theorem}
Again, the expression in~\ref{formulaTH2} is difficult to analyze, but it helps to compute specific (not asymptotic) bounds for some specific values of parameters.
In Table~\ref{tabBP}, we show the values given by the bound for some arbitrary parameters. The bounds are computed the same way as for regular non-bipartite graphs. The bound for bipartite graph is stronger than that of Theorem \ref{BoundReg}. Indeed in order to achieve the same degree, twice as many random matrices are used.

Finally, in order to get theoretical guarantees on biregular graphs, we also prove this upper bound on the merging operation of bipartite regular graphs. Here we consider the non-normalized eigenvalues.

\begin{prop}\label{THbireg}
Let $G$ be a bipartite $d_1$-regular graph, $n$ the size of each partition and $\lambda_2(G)$ the second largest magnitude eigenvalue of $G$. Suppose that $|\lambda_2(G)| \le d_1 \alpha$ for some $ 0 <\alpha < 1$. Let $\gamma$ be an integer that divides $n$ and set $d_2 = \gamma d_1$. Let $H$ be the bipartite $d_1d_2$-biregular graph obtained by merging every $\gamma$ vertices of one of the partitions of $G$ (in an arbitrary order). Let $\lambda_2$ be the second largest eigenvalue of $H$. Then,
$$
|\lambda_2| \le \sqrt{d_1d_2} \alpha.
$$
\end{prop}

\begin{table}[!ht]
\begin{center}
   \begin{tabular}{| l || l | l |}
     \hline
     $(k, d)$   & Th.\ref{BoundBP} & Ramanujan  \\ \hline
     $(14, 30)$   & 0.5787 & 0.3590\\
     $(14, 60)$   & 0.3923 & 0.2560\\
     $(14, 120)$  & 0.2687 & 0.1818\\
     $(20, 40)$   & 0.5741 & 0.3122\\
     $(25, 100)$  & 0.3718 & 0.1989\\
     $(30, 1000)$ & 0.1128 & 0.0632\\
     $(40, 1000)$ & 0.1251 & 0.0632\\
     $(60, 120)$  & 0.5086 & 0.1818\\
     $(200,400)$  & 0.4592 & 0.0998\\
     \hline
   \end{tabular}
\end{center}
\caption{This table shows the upper bound for the expected value of $\mu_2$ computed with Theorem \ref{BoundBP} for a few examples of parameters values
(the dimension $k$ and the degree $d$). Note that the bipartite graphs are constructed by sampling at random $d$ matrices. We display also the Ramanujan's threshold($\frac{2\sqrt{d-1}}{d}$). The first three lines correspond to graphs that we could test in our experiments}
 \label{tabBP}
 \end{table}

\subsection{Proof of Proposition~\ref{WeakBound}}\label{proofP1}

Let $H$ be a group acting transitively on a set $V$ and let $G$ be the Schreier graph of $S \subset H$ acting on $V$ with $S$ a random multiset of $d$ elements. Hence, we obtain an undirected $2d$ regular graph whose size is $n = |V|$. By mapping every element of $V$ to an element of $[\![ 1,\ n ]\!]$ (with some bijective function $f$), we can associate every pair of vertices with a coordinate of a matrix. This way, we can define $M$, the normalized adjacency matrix of $G$. Let 
$$
1 = |\mu_1| \ge |\mu_2| \ge ... \ge |\mu_n|
$$
be its eigenvalues (which are real since the matrix is symmetric). 

Consider a random walk of $2m$ steps starting from vertex $i$. Then the $(i, i)$ coordinate of $M^{2m}$ corresponds to the number of closed walks starting from vertex $i$ of size $2m$ divided by $(2d)^{2m}$, since $(2d)^{2m}$ is the number of paths of size $2m$ starting from $i$. Therefore, this is the probability (denoted $P_{ii}$) of returning to the vertex $i$ after $2m$ steps of the random walk. Since $\textit{Trace}(M^{2m})$ is equal to the sum of all of these quantities and since the expected $(i, i)$ coordinate of $M^{2m}$ is the same for every $i$ we get
\[
\mathbf{E}(\textit{Trace}(M^{2m})) = n \mathbf{E}(P_{11}).
\]
On the other hand, we have 
$$
\sum_{i = 1}^n \mu_i^{2m} = \textit{Trace}(M^{2m}).
$$
Thus, since $\mu_i^{2m} \ge 0$ and $\mu_1 = 1$, we  get $|\mu_2| \le (\textit{Trace}(M^{2m}) - 1)^{1/2m}$, which implies by Jensen's inequality 

\begin{equation}
\label{eq:1}
\mathbf{E}(|\mu_2|) \le (n\mathbf{E}(P_{11}) - 1 )^{1/2m}
\end{equation}

Given the starting vertex (say $v_1 = f^{-1}(1)$), the random walk can be seen as the product of $2m$ group elements that belong to $S$ or their inverses (that is, $S \cup S^{-1}$) chosen uniformly and independently at random. We denote this product of group elements $\omega = w_{2m}.w_{2m - 1} \dots w_2.w_1$ and each $w_i$ represents the choice of a particular neighbor for all vertices. Hence, if the first vertex on the path  is $v_1$, the second  will be $v_2 = s_{i_1}.v_1$ (where $s_{i_1}$ is the value of $w_1$), the third $v_3 = s_{i_2}.v_2$ (where $s_{i_2}$ is the value of $w_2$) and so on; the last vertex of the path is then $\omega.v_1$.

In what follows we call by \emph{literals}  the $2d$ elements of the set 
\[
\{s_1, s^{-1}_1, \dots, s_d, s_d^{-1} \}
\]
where each \emph{letter} $s_j$ (for $j=1,\ldots,d$) appears in two ways: as the literal $s_j$ and as the literal $s_j^{-1}$.
Thus,  a product of elements of $S\cup S^{-1}$ $\omega = w_{2m} w_{2m - 1} \dots w_2 w_1$ can be identified with a sequence of $2m$ literals.

We reuse here the main conceptual idea of the proof in \cite{broder-shamir}. 
The value of $\omega.v_1$ depends on two types of random choices: on the random choice of the word $\omega = w_{2m}.w_{2m - 1} \dots w_2.w_1$ where each literal is chosen at random in  $\{s_1, s^{-1}_1, \dots, s_d, s_d^{-1}\}$, and the random choice of an element in $H$  for each $s_i$. These two choices are independent. It is useful here to sample at first the words $\omega$ and only then choose the elements $s_j$.

We need the following simple lemma in order to estimate the probability of a closed walk. It applies to any transitive group action; this statement also applies if $s_j$ are chosen among  Toeplitz matrices (even though the set of all non-singular Toeplitz matrices with matrix multiplication is not a group). 
\begin{lemma}\label{lemmaUDS}
Let $s$ be an element of a group acting transitively on a set $V$. If letter $s$ appears in a word $\omega$ only once, at some position $i$, then $v_{i+1} = s.v_i$ is chosen uniformly at random among all elements of $V$.
Moreover, for such a word $\omega$, the probability of the event $\omega.v_1 = v_1$ is equal to $1/n$.

The same property is true if $s$ is chosen as a random Toeplitz matrix acting on the set $V$ of non-zero vectors.
\end{lemma}
\begin{proof}
Indeed, if letter $s$ appears in $\omega$ exactly once,  we can rewrite $\omega$ as $AsB$, where  $A$ and $B$ are elements of the group chosen at random (not necessarily uniformly and independently of each other).   The key observation is that $A$ and $B$ are independent from $s$.  So we can at first choose the values of $A$ and $B$. At this stage we still have no information about $s$. 

Further, since the action on $V$ is transitive, for every $v_i$ and a randomly chosen $s$, the values $s.v_i$ are uniformly distributed on $V$.
Thus, as we have fixed  $v_i = B.v_1$, we have $\mathbf{P}(\omega.v_1=v_1) = \mathbf{P}(s.v_i = A^{-1}v_1) = \frac{1}{n}$.

The same is true if $s$ is a Toeplitz matrix. Indeed, a Toeplitz matrix has at least one completely free column. Hence the result of the product $s.v_i$ is uniformly distributed (see Lemma \ref{lemmaUAR} in the next section).
  \qed \end{proof}

The probability that we have just found is computed for a fixed word $\omega$ under the condition that this word contains a letter that appears there exactly once. We denote this condition by $X_1$ (it is an event on $\omega$'s probability space). 
It remains to compute the fraction of words having this property.
Let us remind that $\omega$ can be seen as a word whose symbols (literals) are  taken at random from $\{s_1, s^{-1}_1, \dots s_d, s^{-1}_d\}$.

In order to finish the proof of Proposition \ref{WeakBound}, we bound the number of words in which no letter appears exactly once using an argument similar to that of \cite{alon-roichman}. We observe that for a word that belongs to $\overline{X_1}$ there are at most $m$ different indices $j$ such that $s_j$ or $s_j^{-1}$ (or both) appear in the word. Hence, we have at most $\binom{d}{m}$ ways of choosing those letters in the alphabet. When choosing each letter at random, the probability that all of them are in the right set is $(\frac{m}{d})^{2m}$. Thus,
$$
\mathbf{P}(\overline{X_1}) \le \binom{d}{m}\bigg(\frac{m}{d}\bigg)^{2m} \le \bigg(e\frac{d}{m}\bigg)^{m}\bigg(\frac{m}{d}\bigg)^{2m} = \bigg(e\frac{m}{d}\bigg)^{m}.
$$
The probability we are looking for is then
$$
\mathbf{P}(\omega.v_1=v_1) = \mathbf{P}(\omega.v_1=v_1 | X_1)\mathbf{P}(X_1) + \mathbf{P}(\omega.v_1=v_1 | \overline{X_1})\mathbf{P}(\overline{X_1}) 
$$

$$
\le \mathbf{P}(\omega.v_1=v_1 | X_1) \cdot 1 + 1\cdot \mathbf{P}(\overline{X_1}) \le \frac{1}{n} + \bigg(e\frac{m}{d}\bigg)^{m}.
$$
We set $m = \ln n$ to minimize the bound. Substituting the above quantity back into equation \eqref{eq:1} completes the proof of Proposition \ref{WeakBound}.

\subsection{Proof of Theorem \ref{BoundReg}}\label{proofTh1}

We now turn to the proof of Theorem \ref{BoundReg}.
Let $k$ be a natural integer and $G$ the Schreier graph of $S \subset \textit{GL}_k(\mathbb{F}_2)$ acting on $(\mathbb{F}_2^k)^*$ by matrix-vector product, with $S$ a random multiset of $d$ elements. 
$G$ is thus a $2d$-regular graph of $2^k-1$ vertices.
Let $M$ be the normalized adjacency matrix of $G$.

We reuse the ideas from the preceding section. Let $\omega = w_{2m} \dots w_1$ be a product of matrices chosen at random from $S \cup S^{-1}$. As before, $\omega$ represents a random walk in $G$. Here again, the choice of $S$ and the choice of $\omega$ are independent. Hence we can sample at first the words $\omega$ and only then choose the matrices that are in $S$. 

Following the approach in \cite{broder-shamir}, we prefer not to sample the entire value of each $s_j$ in ``one shot'' but reveal the values of these matrices (better to say, the values of the linear operators corresponding to these matrices) little by little, as it is needed. Thus,  starting at vertex $v_1$, instead of choosing at random in $\textit{GL}_k(\mathbb{F}_2)$ the entire matrix $w_1 = s_{i_1}$, we only determine the result of the product $v_2 = s_{i_1}.v_1$. This choice does not determine completely the matrix $s_i$ but imposes a linear constraint on the  matrix elements of $s_{i_1}$. The same letter $w_1$ may appear in the word $\omega $ several times.  Each time the same letter $w_1$ appears in the word $\omega$ and, therefore, the matrix $s_{i_1}$ is encountered on the path, we must define the action of this matrix on some new vector $x$. We choose the result of $s_{i_1}.x$ by extending the partial definition of $s_{i_1}$, which means an extension of the linear constraints on $s_i$ fixed earlier.
In a similar way, we define step by step the other matrices $s_j$ that are involved in $\omega$.
We need to understand the distribution of the vector  $v_{2m} = \omega.v_1$ that we obtain at the end of this procedure (and the probability of the event $v_{2m} = v_1$).
In the next paragraphs we analyze this distribution. This will lead to a stronger version of Lemma \ref{lemmaUDS} in the case of Schreier graphs of $\textit{GL}_k(\mathbb{F}_2)$.

Consider a matrix $s \in S$ that has been already encountered on the path defined by $\omega$, and we have already defined the action of $s$ on $t$ different vertices. Assume that we encounter the same matrix $s$ once again, and we must  define the product $s.x$ for  some one more vector  $x \in (\mathbb{F}_2^k)^*$. In the permutation model, as stated in \cite{broder-shamir} this would be a uniform distribution over the $n - t$ vertices that have not been earlier assigned to the partially defined permutation $s$. However, in our construction, even if $x$ is totally new to $s$, the result of $s.x$ may not be necessarily undetermined. Indeed, if $x$ is linearly dependent from the vectors that we have already met, we would have
$$
x = \sum_{i = 1}^t \alpha_i x_i,
$$
which is a sum of vectors whose result, when multiplied by $s$, is already known. Thus, 
$$
s.x = \sum_{i = 1}^t \alpha_i s.x_i
$$
would be completely determined by our previous random choices, and would not give any new information about $s$. Intuitively, we would say that the step that leads from $x$ to $s.x$ is not free. In order to characterize formally what it means for a step to be free, we need to introduce the following set: let $s$ be a matrix of $S$, $\omega = w_{2m}\dots w_1$, $v_1$ the starting vertex and $v_{j+1} = w_jv_j$. Then we define
$$\Sigma_s(i) = \Span(\{v_j : j < i, w_j = s\} \cup \{v_{j+1} : j < i, w_j = s^{-1}\}).$$
This is the set of vectors on which the action of $s$ is determined at step $i$. The image set is thus 
$$s.\Sigma_s(i) = \Span(\{v_{j+1} : j < i, w_j = s\} \cup \{v_j : j < i, w_j = s^{-1}\}).$$

This leads to the analogous definition that is presented in \cite{broder-shamir}. 

\begin{definition}[free and forced step]
We consider the $i$-th step in the path. Let $s=w_i$. We say that step $i$ is \emph{forced} when $v_i \in \Sigma_s(i)$. In the opposite case, we say that the step $i$ is \emph{free}.
\end{definition}

Alternatively, instead of saying that a step $i$ is free, we will say that the vector obtained after this step is free (namely the $(i+1)$-th vector, $w_i.v_i$). The following lemma justifies this terminology, and will be used systematically  in the rest of the paper. We prove this lemma for invertible matrices over finite fields of prime size $q$. However, our result uses this only for fields of size $2$.

\begin{lemma}\label{lemmaUAR}
Let $s = w_i$ for a step $i$ and $t$ be the dimension of $\Sigma_s(i)$.
Then, if $v_i \notin \Sigma_s,$ then $ v_{i+1}=s.v_i$ can be chosen uniformly at random among the $q^k - q^t$  vectors that do not belong to $s.\Sigma_s(i)$.
\end{lemma}

\begin{proof}
The choice of a non-degenerate matrix $s$  of size $k\times k$ is the same as the choice of a bijective linear operator from $\mathbb{F}_q^k$ to $\mathbb{F}_q^k$.  To specify a linear operator, we only need to define it on vectors of any basis in  $\mathbb{F}_q^k$.  Let  $x_1, \dots, x_t$ be a basis of $\Sigma_s(i)$. By the assumption of the lemma, vector $v_i$ is linearly independent with $x_1, \dots, x_t$. Therefore, we can let $x_{t+1} = v_i$ and then extend  $x_1, \dots, x_t, x_{t+1}$ to a basis in the space $\mathbb{F}_q^k$ with some $x_{t+2}, \ldots, x_k$.

To define $s$, we should specify one by one linearly independent vectors $y_1=s.x_1, y_2=s.x_2, \ldots, y_k=s.x_k$. We have $q^k-1$ possibilities to choose $y_1$ (any non-zero vector; it is also true for Toeplitz matrices since at least one column is free, hence we can choose $y_1$ as we want), $q^k-q$ possibilities to choose $y_2$ (any vector linearly independent with the fixed $y_1$),  $q^k-q^2$ possibilities to choose $y_3$ (any vector linearly independent with $y_1$ and $y_2$), and so on. In particular, if we have fixed the values $y_i = s.x_i$  
for $i=1,\ldots, t$, then it remains $q^k-q^t$ available options to choose $y_{t+1}$ (which is the same as $v_{i+1}$ in our notation).
  \qed \end{proof}

\paragraph{Remark.} This implies a sort of transitivity property of the group action which is stronger than the simple transitivity, but weaker than the $k$-transitivity: for all $t \le k$, if $(x_1, \dots, x_t)$ and $(y_1, \dots, y_t)$ are two families of linearly independent elements of $\mathbb{F}_q^k$ then, there exists an element $s$ of $\textit{GL}_k(\mathbb{F}_q)$ such that for all $i \le t$, $s.x_i = y_i$.

In order to estimate the total probability of having a closed walk, we subdivide the space of such words in a few subsets (events) and then estimate probabilities of each of them. We chose these subsets so that it will be easier to estimate the conditional probability to get closed path, as we show later. Let us define our events:

\begin{itemize}
\item $X_1$ : ``at least one letter appears in the word exactly once''
\item $X_2$ : $\overline{X_1}$ $\wedge$ ``at least one letter appears exactly twice with same sign''
\item $X_3$ : ``no letter appears once or twice, at least one letter appears exactly three times''
\item $X_4$ : ``no letter appears once, twice, nor three times''
\item $X_2'$ : $\overline{X_1} \wedge \overline{X_3} \wedge \overline{X_4} \wedge$ ``all letters that appear exactly twice have different sign'' \label{SETS}
\end{itemize}

These events form a partition of the set of all possible words, whose size is $(2d)^{2m}$. 
In the proof of Proposition \ref{WeakBound}, we only subdivided the space in two parts: $X_1$ and $\overline{X_1}$. In order to get a tighter bound, we need a more careful analysis of the combinatorial structure of $\omega$ and the corresponding conditional probabilities. We consider again all of our five events. We are going to represent the probability $P(\omega . v_1 = v_1)$ as the sum
\begin{eqnarray}
P(\omega . v_1 = v_1 | \omega\in X_1) \cdot P(X_1) + \ldots +
P(\omega . v_1 = v_1 | \omega\in X_4) \cdot P(X_4) \nonumber \\
{}+ P(\omega . v_1 = v_1 | \omega\in X_2') \cdot P(X_2'). \nonumber
\end{eqnarray}
For $i=1,2,3$ and 4 we estimate separately $P(X_i)$ and  $P(\omega . v_1 = v_1 | \omega\in X_i)$, and we estimate the value of the product $P(\omega . v_1 = v_1 | \omega\in X_2') \cdot P(X_2')$ as a whole.
The sum of these bounds result in the proof of Theorem \ref{BoundReg}.

The events $X_2'$ and $X_4$ together involve one particular event that implies a closed walk with probability $1$. This event  is the ``collapse'' of the whole word to the identity matrix. It happens when iterating the reduction operation ($Ass^{-1}B \mapsto AB$ for all invertible matrices $A$, $s$ and $B$) ends up with the identity matrix. The probability of this event denoted $C$ is analyzed in \cite{broder-shamir} (lemma 2): 
$$
\mathbf{P}(C) = \binom{2m+1}{m}\frac{(2d)^m}{2m+1}\bigg(\frac{1}{2d}\bigg)^{2m} \le \bigg(\frac{2}{d}\bigg)^m.
$$
This is proven by counting the number of well parenthesized words of size $2m$ (Catalan number) with $d$ different type of parenthesis. We wish to bound the probability of having a closed walk when $\omega$'s structure is such that this event cannot happen.

We express the size of these sets using a more general recursive formula.
Let $X_p(c, \ell, d)$ be the size of the set of all words of length $\ell$, on the alphabet that consists of $d$ letters and their negations, such that at least $c$ letters appear (with the positive or negative sign) in this word exactly $p$ times, and the other letters that appear in it have more occurrences. 
Then
\begin{equation}
\label{eq:rec}
X_p(c,\ell, d) = \sum_{i = c}^{\lfloor\frac{\ell}{p}\rfloor} \binom{d}{i}
\prod_{j = 0}^{i-1}2^p\binom{\ell-jp}{p}X_{p+1}(0,\ell-ip, d-i)
\end{equation}
Indeed, $2^p\binom{\ell}{p}$ is the number of ways of placing $p$ times the same letter in a word of size $\ell$ (each letter can have positive or negative sign).
Thus, $\prod_{j = 0}^{i-1}2^p\binom{\ell-jp}{p}$ is the number of ways of repeating $i$ times this operation while removing at every step $p$ free places. It simplifies as follows 
$$
\prod_{j = 0}^{i-1}2^p\binom{\ell-jp}{p} = \bigg(\frac{2^p}{p!}\bigg)^i\prod_{j = 0}^{i-1}\frac{(\ell-jp)!}{(\ell- jp - p)!} = \bigg(\frac{2^p}{p!}\bigg)^i\frac{\ell!}{(\ell-ip)!}
$$

Then, once the $i$ different letters are placed $p$ times (we have $\binom{d}{i}$ ways of choosing the $i$ letters), we know that the $d-i$ other different letters will appear either 0 or more than $p+1$ times, and the remaining spaces will be $\ell - ip$. This explains the recursive call in \ref{eq:rec}. Note that $X_p(0,0,d) = 1$ because we only have one way of placing no letters in a word of size 0. Moreover, if $p > \ell$ then $X_p(c,\ell,d) = 0$, since the $p$ letters cannot fit in the word. Using these observations, we have a complete recursive definition of $X_p(c,\ell, d)$,  
\[
X_p(c,\ell, d) = \left\{
    \begin{array}{ll}
        1 \mbox{ if } c = 0 \mbox{ and } \ell = 0; \\
        0 \mbox{ if } p > \ell; \\
        \sum_{i = c}^{\lfloor\frac{\ell}{p}\rfloor} \binom{d}{i}
\bigg(\frac{2^p}{p!}\bigg)^i\frac{\ell!}{(\ell-ip)!}X_{p+1}(0,\ell-ip, d-i) \mbox{ otherwise.}
    \end{array} 
\right.
\]
Then we get
$$
|X_1| = X_1(1, 2m, d),
$$
$$
|X_3| = X_3(1, 2m, d)
$$
and
$$
|X_4| = X_4(0, 2m, d).
$$
One can notice that $|X_3\sqcup X_4| = |X_3|+ |X_4| = X_3(0, 2m, d)$.

Since there are $4^i - 2^i$ possibilities for choosing the sign of $i$ pairs so that the letters of at least one of these pairs have same sign (there are $2^i$ ways of not verifying that), the number of ways of placing these pairs in the word is 
\[
\binom{d}{i}(4^i - 2^i)\prod_{j=0}^{i-1}\binom{2m-2j}{2} = \binom{d}{i}(2^i - 1)\frac{(2m)!}{(2m-2i)!}.
\]
Hence
\[
|X_2| = \sum_{i = 1}^m \binom{d}{i} (2^i - 1)\frac{(2m)!}{(2m-2i)!} X_3(0, 2m-2i, d-i).
\]
Similarly, there are $2^i$ ways of choosing the sign of $i$ pairs of letter so that all pairs are of different sign. Thus we get
\[
|X_2'| = \sum_{i = 1}^m \binom{d}{i} \frac{(2m)!}{(2m-2i)!} X_3(0,2m-2i, d-i).
\]
This can be summarized with the relation $X_2(1, 2m, d) = |X_2| + |X_2'|$.

We have already explained that $\mathbf{P}(\omega.v_1 = v_1|X_1) = \frac{1}{n}$ (see Lemma \ref{lemmaUDS}). It remains to bound this probability conditioned to the other events. We start with $X_2$. 

\begin{lemma}\label{lemmaX2}
$$
\mathbf{P}(\omega.v_1 = v_1 | X_2) \le \frac{2}{n}.
$$
\end{lemma}
\begin{proof}
Let $s$ be the matrix that appears twice with same sign. The word is then of the form $AsBsC.v_1 = v_1$ with $A, B$ and $C$ some invertible matrices of known coefficients.
We can rewrite this equation as $sBs.x = y$ with $x$ and $y$ two determined vectors ($x=C.v_1$ and $y = A^{-1}.v_1$). 
It is useful to name the different vectors of the product:  
$$
\underbrace{s\ \overbrace{B\ \underbrace{s.x}_{y'}}^{x'}}_{y''} = y.
$$
Since we are in the field of size two there is no non-trivial pairs of parallel vectors. Hence the step that leads to $y''$ is free only if $x' \neq x$. In a larger field ($q > 2$), for $y''$ to be free, it is necessary that $x \neq \alpha x'$ for all non-zero $\alpha \in \mathbb{F}_q$. By taking $q=2$, a lot of case-by-case analysis is avoided.

Because $y'$ is necessarily free and since $x$ and $x'=B.y'$ are independent, $\mathbf{P}(x' = x) = \frac{1}{n}$. Then, if $x' = x$, we have $y'' = y'$. This is the probability that $y''$ is forced. If this is not the case, namely if $x' \neq x$ (which happens with probability $\frac{n-1}{n}$), then the probability for $y''$ to be equal to $y$ is at most $\frac{1}{n-1}$ ($y''$ cannot be equal to $y'$ since both steps are free). Therefore, 
$$
\mathbf{P}(\omega.v_1 = v_1 | X_2) \le \frac{1}{n} + \frac{n-1}{n}\frac{1}{n-1} = \frac{2}{n}. 
$$
  \qed \end{proof}

Now we bound the probability $\mathbf{P}(\omega v_1 = v_1 | X_3)$. We proceed the same way as above, by distinguishing the cases where the final step is free or not. We prove the following claim:
\begin{lemma}\label{lemmaX3}
$$
\mathbf{P}(\omega.v_1 = v_1 | ``at\ least\ a\ letter\ appears\ exactly\ three\ times") \le \frac{5}{n}.
$$

In particular, we have
$$
\mathbf{P}(\omega.v_1 = v_1 | X_3) \le \frac{5}{n}.
$$
\end{lemma}
\begin{proof}
Under the $X_3$ condition, the word can take four forms that will be analyzed separately:
\begin{itemize}
\item $sBsCs.x = y$
\item $s^{-1}BsCs.x = y$
\item $sBs^{-1}Cs.x = y$
\item $sBsCs^{-1}.x = y$
\end{itemize}
Any other form can be turned into one of the above by switching $s$ with $s^{-1}$, which does not change the argument.
The different possibilities can be summarized by writing $s_1Bs_2Cs_3.x = y$; at most one of $s_1, s_2, s_3$ is $s^{-1}$ and the others are $s$. We will use the notation below to treat all four cases:
$$
\underbrace{s_1\ \overbrace{B\ \underbrace{s_2\ \overbrace{C\ \underbrace{s_3.x}_{y'}}^{x'}}_{y''}}^{x''}}_{y'''} = y.
$$

We start with the case in which there is no $s^{-1}$ in $\omega$. Here, the step that leads to $y'''$ is free only if $x'' \neq x$, $x'' \neq x'$ and $x'' \neq x' + x''$ (that is, $x''$ is not a linear combination of $x'$ and $x$).
Since $x'$ and $x$ are independent,  $\mathbf{P}(x'=x) =  \frac{1}{n}$. Hence, with probability $\frac{n-1}{n}$ we get that $y''$ is free, which means that $x'' = B.y''$ is uniformly distributed among the $n-1$ vectors different from $B.y'$. There are three values for $x''$ that make the final step forced and they are equally likely, thus $\mathbf{P}(y'''\text{ is forced }| x' \neq x) \le \frac{3}{n-1}$. The opposite case happens with probability $\frac{n-4}{n-1}$. Then $\mathbf{P}(y''' = y) \le \frac{1}{n-3}$. To illustrate the reasoning, we can represent all possible cases as a tree:

\begin{center}
\begin{tikzpicture}[roundnode/.style={draw,  minimum width=2pt, inner sep=1pt}]
        \node[roundnode](1){};
        \node[roundnode](2)[left of=1,yshift=-1.5cm,xshift=-0.4cm]{$x'=x$};
        \node[roundnode](3)[right of=1,yshift=-1.5cm,xshift=0.4cm]{$x'\neq x$};
        \node[roundnode](4)[left of=3, yshift=-1.5cm,xshift=-0.4cm]{$x''=x$};
        \node[roundnode](5)[right of=3, yshift=-1.5cm,xshift=0.4cm]{$x''\neq x$};
        \node[roundnode](6)[left of=5, yshift=-1.5cm,xshift=-0.4cm]{$x''=x'$};
        \node[roundnode](7)[right of=5, yshift=-1.5cm,xshift=0.4cm]{$x''\neq x'$};
        \node[roundnode](8)[left of=7, yshift=-1.5cm,xshift=-0.4cm]{$x''=x+x'$};
        \node[roundnode](9)[right of=7, yshift=-1.5cm,xshift=0.4cm]{$x''\neq x+x'$};

        \draw [-](1) -- (2) node [midway, fill=white]{$\frac{1}{n}$};
        \draw [-](1) -- (3) node [midway, fill=white]{$\frac{n-1}{n}$};
        \draw [-](3) -- (4) node [midway, fill=white]{$\frac{1}{n-1}$};
        \draw [-](3) -- (5) node [midway, fill=white]{$\frac{n-2}{n-1}$};
        \draw [-](5) -- (6) node [midway, fill=white]{$\frac{1}{n-2}$};
        \draw [-](5) -- (7) node [midway, fill=white]{$\frac{n-3}{n-2}$};
        \draw [-](7) -- (8) node [midway, fill=white]{$\frac{1}{n-3}$};
        \draw [-](7) -- (9) node [midway, fill=white]{$\frac{n-4}{n-3}$};

\end{tikzpicture}
\end{center}
The rightmost leaf corresponds to $y'''$ being free which gives a probability $\frac{1}{n-3}$ of having a closed walk. Therefore, 
$$
\mathbf{P}(sBsCs.x=y) \le \frac{1}{n} + \frac{n-1}{n}\bigg(\frac{3}{n-1} + \frac{n-4}{n-1}\frac{1}{n-3}\bigg) \le \frac{5}{n}. 
$$

Now, consider $s_3 = s^{-1}$. Then $y''$ is forced if $x' = y'$, but those two vectors are correlated, so we cannot bound the probability of this event. We will consider both cases and take the probability of the most likely event as a bound.
If $y' = x'$, then $y''$ is forced, which implies that $y'' = x$. In this case, if $x'' = y'$ we have $y''' = x$. However, $x'' = B.x$, which is independent from $y'$ (which is from a free step). Hence, the probability for them to be equal is $\frac{1}{n}.$ In the opposite case, $y'''$ is free, which gives a total probability of this branch of $\frac{2}{n}$.
We now suppose that $y''$ is free. Then, with probability $\frac{3}{n-1}$, $y'''$ is forced. In the other case, $y'''$ is equal to $y$ with probability at most $\frac{1}{n-3}$. Here is the tree of all possible cases:

\begin{center}
\begin{tikzpicture}[roundnode/.style={draw,  minimum width=2pt, inner sep=1pt}]
        \node[roundnode](1){};
        \node[roundnode](2)[left   of=1,yshift=-1.5cm,xshift=-1.2cm]{$x'=y'$};
        \node[roundnode](10)[right of=2,yshift=-1.5cm,xshift=0cm]{$x''\neq y'$};
        \node[roundnode](11)[left  of=2, yshift=-1.5cm,xshift=-0cm]{$x''=y'$};
        \node[roundnode](3)[right  of=1,yshift=-1.5cm,xshift=1.2cm]{$x'\neq y'$};
        \node[roundnode](4)[left   of=3, yshift=-1.5cm,xshift=-0.4cm]{$x''=x'$};
        \node[roundnode](5)[right of=3, yshift=-1.5cm,xshift=0.4cm]{$x''\neq x'$};
        \node[roundnode](6)[left  of=5, yshift=-1.5cm,xshift=-0.4cm]{$x''=y'$};
        \node[roundnode](7)[right of=5, yshift=-1.5cm,xshift=0.4cm]{$x''\neq y'$};
        \node[roundnode](8)[left  of=7, yshift=-1.5cm,xshift=-0.4cm]{$x''=y'+x'$};
        \node[roundnode](9)[right of=7, yshift=-1.5cm,xshift=0.4cm]{$x''\neq y'+x'$};
        
        \draw [-](1) -- (2) node [midway]{};
        \draw [-](1) -- (3) node [midway]{};
        \draw [-](2) -- (10) node [midway, fill=white]{$\frac{n-1}{n}$};
        \draw [-](2) -- (11) node [midway, fill=white]{$\frac{1}{n}$};
        \draw [-](3) -- (4) node [midway, fill=white]{$\frac{1}{n-1}$};
        \draw [-](3) -- (5) node [midway, fill=white]{$\frac{n-2}{n-1}$};
        \draw [-](5) -- (6) node [midway, fill=white]{$\frac{1}{n-2}$};
        \draw [-](5) -- (7) node [midway, fill=white]{$\frac{n-3}{n-2}$};
        \draw [-](7) -- (8) node [midway, fill=white]{$\frac{1}{n-3}$};
        \draw [-](7) -- (9) node [midway, fill=white]{$\frac{n-4}{n-3}$};

\end{tikzpicture}
\end{center}
The branches whose probabilities are close to 1 corresponds to the cases when $y'''$ is free. Thus we have 
$$
\mathbf{P}(sBsCs^{-1}.x=y) \le \max\bigg(\frac{2}{n}, \frac{3}{n-1} + \frac{n-4}{n-1}\frac{1}{n-3}\bigg) \le \frac{4}{n-1}. 
$$

When $s_2 = s^{-1}$ the choice of $y''$ is forced if $x' = y'$, which implies $x'' = B.x$. Then $s$ is defined only on $x$. If $B.x = x$ then, since $y'$ is free, $\mathbf{P}(y''' = y) = \mathbf{P}(y' = y) = \frac{1}{n}$. Otherwise, $y'''$ is free, and therefore $\mathbf{P}(y''' = y) = \frac{1}{n-1}$.
On the other hand, if $y''$ is free, then $s$ is defined on $x$ and $y''$. Since $y''$ is free, $x''$ and $x$ are independent, thus $\mathbf{P}(x'' = x) = \frac{1}{n-1}$. Moreover, if $x'' = y''$,  $\mathbf{P}(y'''=y)=\mathbf{P}(x'=y)=\frac{1}{n}$. If $x'' = y''+x$, we have $y'''=x'+y' = (C+Id_k)y'$ which is independent from $y$. Hence, $\mathbf{P}(y'''=y)=\frac{1}{n}$. Otherwise, since three vectors are excluded, $y'''$ is free with probability at most $\frac{n-3}{n}$. If so, $\mathbf{P}(y''=y)=\frac{1}{n-3}$. As before, the case study can be illustrated with a tree:

\begin{center}
\begin{tikzpicture}[roundnode/.style={draw,  minimum width=2pt, inner sep=1pt}]
        \node[roundnode](1){};
        \node[roundnode](2)[left   of=1,yshift=-1.5cm,xshift=-1.2cm]{$x'=y'$};
        \node[roundnode](10)[right of=2,yshift=-1.5cm,xshift=0cm]{$x''\neq x$};
        \node[roundnode](11)[left  of=2, yshift=-1.5cm,xshift=-0cm]{$x''=x$};
        \node[roundnode](3)[right  of=1,yshift=-1.5cm,xshift=1.2cm]{$x'\neq y'$};
        \node[roundnode](4)[left   of=3, yshift=-1.5cm,xshift=-0.4cm]{$x''=x$};
        \node[roundnode](5)[right of=3, yshift=-1.5cm,xshift=0.4cm]{$x''\neq x$};
        \node[roundnode](6)[left  of=5, yshift=-1.5cm,xshift=-0.4cm]{$x''=y''$};
        \node[roundnode](7)[right of=5, yshift=-1.5cm,xshift=0.4cm]{$x''\neq y''$};
        \node[roundnode](8)[left  of=7, yshift=-1.5cm,xshift=-0.4cm]{$x''=y''+x$};
        \node[roundnode](9)[right of=7, yshift=-1.5cm,xshift=0.4cm]{$x''\neq y''+x$};
        
        \draw [-](1) -- (2) node [midway]{};
        \draw [-](1) -- (3) node [midway]{};
        \draw [-](2) -- (10) node [midway]{};
        \draw [-](2) -- (11) node [midway]{};
        \draw [-](3) -- (4) node [midway, fill=white]{$\frac{1}{n-1}$};
        \draw [-](3) -- (5) node [midway, fill=white]{$\frac{n-2}{n-1}$};
        \draw [-](5) -- (6) node [midway]{};
        \draw [-](5) -- (7) node [midway]{};
        \draw [-](7) -- (8) node [midway, fill=white]{$\frac{1}{n-2}$};
        \draw [-](7) -- (9) node [midway, fill=white]{$\frac{n-3}{n-2}$};

\end{tikzpicture}
\end{center}

Thus, we get
$$
\mathbf{P}(sBs^{-1}Cs.x=y) \le  \max \bigg(\frac{2}{n-1},  \frac{3}{n-1}+\frac{n-2}{n-1}\bigg(\frac{n-3}{n-2} \frac{1}{n-3}\bigg)\bigg) = \frac{4}{n-1}. 
$$

Lastly, we consider the case $s_1 = s^{-1}$. Here, $y''$ is forced when $x'= x$. Those are not correlated, so this event happens with the probability $\frac{1}{n}$. In the other case, $s^{-1}$ is defined on $y''$ and $y'$, which are random. If $x'' = y''$, we have $y''' = x'$ which is equal to $y$ with probability less than $\frac{1}{n-1}$. Since $y''$ is free, $y'$ and $y''$ are independent, hence $\mathbf{P}(x'' = y') = \mathbf{P}(x'' = y' + y'') = \frac{1}{n-1}$. If $y'''$ is free, it can take any value with probability $\frac{1}{n-3}$.
The last tree of possible cases is then

\begin{center}
\begin{tikzpicture}[roundnode/.style={draw,  minimum width=2pt, inner sep=1pt}]
        \node[roundnode](1){};
        \node[roundnode](2)[left of=1,yshift=-1.5cm,xshift=-0.4cm]{$x'=x$};
        \node[roundnode](3)[right of=1,yshift=-1.5cm,xshift=0.4cm]{$x'\neq x$};
        \node[roundnode](4)[left of=3, yshift=-1.5cm,xshift=-0.4cm]{$x''=y''$};
        \node[roundnode](5)[right of=3, yshift=-1.5cm,xshift=0.4cm]{$x''\neq y''$};
        \node[roundnode](6)[left of=5, yshift=-1.5cm,xshift=-0.4cm]{$x''=y'$};
        \node[roundnode](7)[right of=5, yshift=-1.5cm,xshift=0.4cm]{$x''\neq y'$};
        \node[roundnode](8)[left of=7, yshift=-1.5cm,xshift=-0.4cm]{$x''=y'+y''$};
        \node[roundnode](9)[right of=7, yshift=-1.5cm,xshift=0.4cm]{$x''\neq y'+y''$};

        \draw [-](1) -- (2) node [midway, fill=white]{$\frac{1}{n}$};
        \draw [-](1) -- (3) node [midway, fill=white]{$\frac{n-1}{n}$};
        \draw [-](3) -- (5) node [midway]{};
        \draw [-](3) -- (4) node [midway]{};
        \draw [-](5) -- (6) node [midway, fill=white]{$\frac{1}{n-2}$};
        \draw [-](5) -- (7) node [midway, fill=white]{$\frac{n-3}{n-2}$};
        \draw [-](7) -- (8) node [midway, fill=white]{$\frac{1}{n-3}$};
        \draw [-](7) -- (9) node [midway, fill=white]{$\frac{n-4}{n-3}$};

\end{tikzpicture}
\end{center}
Hence we have
$$
\mathbf{P}(s^{-1}BsCs.x = 1) \le \frac{1}{n} +  \frac{n-1}{n}\max\bigg(\frac{1}{n-2}, \frac{1}{n-2}+\frac{n-4}{n-2}\frac{1}{n-3}\bigg) \le \frac{5}{n}. 
$$
By taking the maximum of all these bounds, we conclude the proof.
  \qed \end{proof}

It remains to bound $\mathbf{P}(\omega.v_1 = v_1 | X_2')\mathbf{P}(X_2')$. To simplify the notations we set $x_3(i) = X_3(1, 2m-2i, d-i)$ and $x_4(i) = X_4(0, 2m-2i, d-i)$.
We need to prove the following  statement.
\begin{lemma}\label{lemmaX2'}
\begin{multline}
\mathbf{P}(\omega.v_1 = v_1 | X_2')\mathbf{P}(X_2') \le
\\ \bigg(\frac{1}{2d}\bigg)^{2m} \sum_{i = 1}^m \binom{d}{i} \frac{(2m)!}{(2m-2i)!} \bigg[\frac{x_3(i) + x_4(i)}{n}  + \frac{2^i}{(i+1)!}\bigg(\frac{5}{n}x_3(i) + x_4(i)\bigg)\bigg].
\end{multline}

\end{lemma}
\begin{proof}

Before we proceed with the proof of this lemma we stress again that in this statement we do not bound separately $\mathbf{P}(X_2')$ and $\mathbf{P}(\omega.v_1 = v_1 | X_2')$, we estimate directly the product of these two probabilities, which equals to the probability of the event
$$
\mathbf{P}(\omega.v_1 = v_1  \text{ and } \omega\in X_2').
$$

The probability is taken, as usual, over the random choice of a word $\omega$ of $2m$ letters and the  random  choice of invertible matrices assigned to the letters of this alphabet.

We start the proof with two claims.

\smallskip
\noindent
\emph{Claim 1:
Assume that the word $\omega$ contains letters $t$ and $s$ exactly twice, and each of these letters appears once with the positive and once with the negative sign,
and these letters interleave:
\begin{equation}\label{tsts}
\omega  = \ldots t \ldots s \ldots t^{-1} \ldots s^{-1} \ldots
\end{equation}
Then the probability to get a closed walk corresponding to the path $\omega$ (probability taken over the choice of matrices for each letter in the alphabet)
is equal to $1/n$.
The claim remains true if we swap the positions of the pair of letters $s$ and $s^{-1}$ and/or of the pair of letters $t$ and $t^{-1}$.
}

\smallskip

\begin{proof}[proof of Claim 1]
Words from $X_2'$ are all of the form $AsBs^{-1}C$ with $A$, $B$ and $C$ some invertible matrices. Hence, we wish to estimate the probability of the event $sBs^{-1}.x = y$, with $x = C.v_1$, and $y = A^{-1}.v_1$. We use the notation
$$
\underbrace{s\ \overbrace{B\ \underbrace{s^{-1}.x}_{y'}}^{x'}}_{y''} = y.
$$
Here, the matrix $t$ is a factor of $B$ (hence $B = \dots t \dots$). 
We first suppose that $x=y$. Then, if $x' = y'$ we have $y'' = x = y$. Since $t$ appears in $B$, $x'$ is independent of $y'$, and thus $\mathbf{P}(y'=x')=\frac{1}{n}$ (because this is the first time $t$ is used in the path). In the opposite case, we have $y'' \neq y$, and the path cannot be closed. 
 
Now we suppose $x\not=y$. Then if $y'=x'$ we have $y'' = x \neq y$. If $y' \neq x'$ (which happens with probability $\frac{n-1}{n}$, $y''$ is free, and its value is uniformly distributed among the $n-1$ remaining vectors.
Therefore, when we have this configuration of random matrices in $\omega$, the probability of having a closed walk is $\frac{1}{n}$.
  \qed \end{proof}

It can be noticed that here, the fact that $t$ appears with different sign is not used.

\smallskip
\noindent
\emph{Claim 2:
Let us take the set of $2i$ literals 
\[
\{ s_1, s_1^{-1}, \ldots, s_i, s_i^{-1}\}
\]
and consider the set of all words of length $(2i)$ composed of these literals (each one should be used exactly once). We claim that the fraction of words that represent a well formed structure of $i$ pairs of parentheses, where each pairs is associated with some pair of literals  $(s_j,s^{-1}_j)$ or $(s^{-1}_j,s_j)$, is equal to $\frac{2^i}{(i+1)!}$.
}
\begin{proof}[proof of Claim 2]
In general, we have $(2i)!$ different ways to distribute $(2i)$ literals among $(2i)$ positions. Let us count the fraction of permutations where the literals form a  structure of $i$ pairs of parentheses.
The number of well parenthesized words (with one type of parentheses) of size $2i$ is the Catalan number $\binom{2i+1}{i}\frac{1}{2i+1}$.  We have  $i!$ ways to assign each pair of parentheses with one of $i$ types of literals, and  $2^i$ to chose the signs in each pairs ($\ldots s_j \ldots s_j^{-1} \ldots$ or $\ldots s_j^{-1} \ldots s_j \ldots$ for each of $i$ pairs). Hence, the proportion of the well parenthesized words is 
\[
\frac{\binom{2i+1}{i}\frac{1}{2i+1}i!2^i}{(2i)!}=\frac{(2i+1)!i!2^i}{(2i+1)!i!(i+1)!}=\frac{2^i}{(i+1)!}.
\]

Observe that this computation makes sense even for $i=1$. Indeed, if we are given only one pair of letters $(s, s^{-1})$, then  we have exactly two possible words ($s s^{-1}$ and $s^{-1}s$), and each of them is  well parenthesized). So  the proportion of the of well parenthesized words is  $\frac{2^1}{(1+1)!} = 1$.
  \qed \end{proof}

It is easy to see that the absence of pattern \ref{tsts} is equivalent to having such well formed structure of parenthesis. Words from $X_2'$ are all of the form $AsBs^{-1}C$, where $A,B,C$ are obtained by a product of  several matrices corresponding to other letter (excluding $s$). Observe that we must take into consideration the case when $B$ is  the identity operator (and two copies of $s$ with the opposite signs collapse with each other).

Let us proceed with the proof of the lemma. By definition, in each word $\omega\in X_2'$ all letters that appear exactly twice must have different signs. In what follows we denote $i$ that number of letters that appear in $\omega$ exactly twice.  For a fixed $i$, to specify a word $\omega$ where $i$ letters appear twice (with opposite signs) and the other letters appears at least three times, we should
\begin{itemize}
\item choose $i$ letters among $d$ (those who appear exactly twice), which gives  $\binom{d}{i}$ combinations;
\item choose $2i$ positions in the word $\omega$ of length $2m$ where we place the letters that appear twice, which gives  $\binom{2m}{2i}$ combination;
\item fix a permutations of the $(2i)$ literals on the chosen $(2i)$ positions, which gives $(2i)!$ combinations;
\item fill the remaining $(2m-2i)$ positions of $\omega$ with other letters, using each letter at least three times; we subdivide these combinations into two sub-cases:
\begin{itemize}
\item there is at least one letter that is used \emph{exactly} three times;  we have $x_3(i) = X_3(1, 2m-2i, d-i)$ possibilities to do it;
\item there is no letter that is used exactly three times, i.e., each letter (besides the $i$ letters that were used twice) must be used at least four times;  we have  $x_4(i) = X_4(0, 2m-2i, d-i)$ possibilities to fill in this way the remaining $(2m-2i)$  positions.
\end{itemize}
\end{itemize}
The  $i$ pairs of letters in $\omega$ contain the pattern \eqref{tsts} may contain or not contain the pattern \eqref{tsts}. By Claim~2, the latter is the case for the fraction $\frac{2^i}{(i+1)!}$ of all $\omega$
(with $i$ pairs) and, respectively,  the former is the case  for the fraction $1-\frac{2^i}{(i+1)!}$ of these words.

If the $i$ pairs of letters in $\omega$ contain the pattern \eqref{tsts}, then by Claim~1 the probability that $\omega$ provides a closed path  is at most $\frac1n$ (probability taken over the choice of matrices for each letter in the alphabet). 
Since we have in total $(2d)^{2m}$ words $\omega$, this case contributes to the  resulting probability $\mathbf{P}(\omega.v_1 = v_1  \text{ and } \omega\in X_2')$  at most by
\[ 
\bigg(\frac{1}{2d}\bigg)^{2m} \binom{d}{i}  \binom{2m}{2i}  (2i)! \bigg(1-\frac{2^i}{(i+1)!}\bigg) (x_3(i)+x_4(i)) \cdot \frac1n
\]
(in what follows we bound $1-\frac{2^i}{(i+1)!}$ by $1$).

If the $i$ pairs of letters in $\omega$ do not contain the pattern \eqref{tsts} but one of the other letters appears in $\omega$  exactly three times, then
the probability of having a closed path is at most $\frac{5}{n}$, as shown in  Lemma~\ref{lemmaX3}.  This case contributes to the  resulting probability   at most by
\[
\bigg(\frac{1}{2d}\bigg)^{2m} \binom{d}{i}  \binom{2m}{2i}  (2i)! \cdot  \frac{2^i}{(i+1)!} \cdot x_3(i)\cdot  \frac5n
\]
At last, if  $\omega$ does not contain the pattern \eqref{tsts}  and all other letters appearing in $\omega$ are used more than three times, then we trivially bound the probability to have a closed path
by $1$. This contributes to the  resulting probability  by
\[
\bigg(\frac{1}{2d}\bigg)^{2m} \binom{d}{i}  \binom{2m}{2i}  (2i)! \cdot  \frac{2^i}{(i+1)!} \cdot x_4(i).
\]
Summing these quantities for all possible  values of $i$ and observing that \\$\binom{2m}{2i}  (2i)! = \frac{(2m)!}{(2m-2i)!}$, we obtain the statement of the lemma.
  \qed \end{proof}

We proceed with similar bounds for the other sets of words:

$$
\mathbf{P}(\omega.v_1 = v_1 | X_1)\mathbf{P}(X_1) \le \frac{1}{n}\bigg(\frac{1}{2d}\bigg)^{2m} |X_1|,
$$

$$
\mathbf{P}(\omega.v_1 = v_1 | X_2)\mathbf{P}(X_2) \le \frac{2}{n}
\bigg(\frac{1}{2d}\bigg)^{2m}|X_2|,
$$

$$
\mathbf{P}(\omega.v_1 = v_1 | X_3)\mathbf{P}(X_3) \le \frac{5}{n}
\bigg(\frac{1}{2d}\bigg)^{2m}|X_2|,
$$
and
$$
\mathbf{P}(\omega.v_1 = v_1 | X_4)\mathbf{P}(X_4) \le \bigg(\frac{1}{2d}\bigg)^{2m} |X_4|.
$$
The sum of these expressions is larger than $P_{11}$ defined at the beginning of this section. By replacing it in equation \ref{eq:1}, we complete the proof. It is easy to see that the rough bounding used in the proof of Proposition \ref{WeakBound} gives a larger bound than that of Theorem~\ref{BoundBP}.

\subsection{Proof of Theorem \ref{BoundBP}}

We now can adapt this proof to get a similar bound for $d$-regular bipartite graphs. Let $G$ be a bipartite Schreier graph of $\textit{GL}_k(\mathbb{F}_2)$ acting on $\mathbb{F}_2^k$ with respect to the sub-multiset $D$ and $M$ be its normalized adjacency matrix. In order to associate its coordinates to vertices we can proceed as in the preceding section by mapping surjectively  $[\![ 1,\ 2(2^k-1) ]\!]$ to $(\mathbb{F}_2^k)^*$, taking care of distinguishing the vectors of the first and the second partition. Here, we set $2n = 2(2^k-1)$, the number of vertices in the graph. Let us start by adapting the trace method to the bipartite graphs. One can remark that the adjacency matrix of $G$ is of the form
\[
M = \left(\begin{array}{@{}c|c@{}}
  0 & A \\ \hline \,^t A & 0
\end{array}\right)
\]
where $^t A$ is the transposition of $A$. In a bipartite graph, it is known that the spectrum $|\mu_1| \ge \dots \ge |\mu_{2n}|
$ is symmetric with respect to zero. Hence for $1\le i\le n$, we have $|\mu_{2i+1}| = |\mu_{2(i+1)}|$. This way we get
$$
\sum_{i = 0}^{n-1} 2 \mu_{2i+1}^{2m} = \textit{Trace}(M^{2m}).
$$
In order to study the spectral gap, the relevant quantity to estimate is then $|\mu_3| = |\mu_4|$. Since $|\mu_1| = |\mu_2| = 1$, we thus obtain
$$
2\mu_3^{2m} \le \sum_{i = 1}^{n-1} 2 \mu_{2i+1}^{2m} = \textit{Trace}(M^{2m}) - 2.
$$

As we have seen in section 2.3.4, the expected value of $\textit{Trace}(M^{2m})$ is the sum of the probability of getting a closed path of size $2m$, starting on each vertex. We note this probability $P_{ii}$ for the vertex $i$. It is the same for every vertex, hence we get, by using Jensen's inequality
\begin{equation}
\label{eq:2}
\mathbf{E}(|\mu_3|) \le \bigg(\frac{1}{2}(\mathbf{E}(\textit{Trace}(M^{2m})) - 2)\bigg)^{\frac{1}{2m}} = (nP_{11} - 1)^{\frac{1}{2m}}.
\end{equation}
One can notice that here, $n$ is the size of the partition, not the size of the graph. Indeed, with an even number of steps, the path must end in the same partition as it started, which eliminates half of the vertices.

We first explain why the construction for even degree regular bipartite graphs gives the same bound as Theorem \ref{BoundReg}. Here, a random walk can be represented as a sequence of matrices of $D \cup D^{-1}$. This is because every vertex $x$ is connected to $s.x$ and $s^{-1}.x$. Each element of the sequence is chosen independently of the others. It is then easy to see that the structure of the walk is exactly the same as in the non-bipartite case: a uniformly random sequence of $2m$ matrices from $D \cup D^{-1}$. The same proof can then be applied to this sequence, the elements of the sequence will then behave the same way as in the preceding section.
 
However, some work needs to be done for graphs of odd degree. In order to apply here a similar reasoning as in the previous section, we need to understand what a random walk in $G$ looks like in terms of the matrices of $D$.

\begin{figure}[ht]
\begin{center}
\begin{tikzpicture}[roundnode/.style={circle, draw,  minimum width=2pt, inner sep=1pt}]
        \node[roundnode](1){$x$};
        \node[roundnode](2)[below of=1]{};
        \node[roundnode](3)[below of=2]{};
        \node[roundnode](4)[below of=3]{$z$};
        \node[roundnode](a)[right of=1, xshift=1cm]{};
        \node[roundnode](b)[right of=2, xshift=1cm]{$y$};
        \node[roundnode](c)[right of=3, xshift=1cm]{};
        \node[roundnode](d)[right of=4, xshift=1cm]{};
        
        node at (1, 2) {\vdots};

        \draw [->](1) -- (b) node [midway, fill=white]{$S_i$};

        \draw [->](b) -- (4) node [midway, fill=white] {$S_j^{-1}$};

\end{tikzpicture}
\caption{The steps of the random walk work by pairs of matrices of $D$; the first one brings us on a vertex of the right hand side, the other is for the way back.}
\label{fig:3}
\end{center}
\end{figure}
In a bipartite regular graph obtained by our construction, a random walk of size $2m$ is a sequence $\omega$ of elements of $D$ chosen independently at random. As usual, every edge in the graph corresponds to some invertible matrix. In the case of a bipartite graph we assume that the multiplication by the matrices transforms the vertices in the left part into the vertices of the right part. Thus, in a random walk on such a graph, 
the matrices that appear in an odd position in $\omega$ are taken with positive sign, while the matrices that appear in an even step are taken with negative sign, see Fig.~\ref{fig:3}.

We wish to proceed as in the preceding section, by partitioning the set of possible sequences (words) so that we can analyze the probability of having a closed walk conditioned to the sets this partition. In order to reuse the results above, we choose a similar partitioning. There is a minor difference: now the signs of matrices appearing in $\omega$  are fixed (positive for odd position and negative for the others). Hence, the number of possible words is $d^{2m}$. We define the sets before determining their size:

\begin{itemize}
\item $Y_1$ : ``at least one letter appears exactly once"
\item $Y_2$ : $\overline{Y_1}$ $\wedge$ "at least one letter appears exactly twice"
\item $Y_3$ : ``no letter appears once or twice, at least one letter appears exactly three times"
\item $Y_4$ : ``no letter appears once, twice, nor three times"
\end{itemize}

Up to this point, it is not hard to understand why the bound of Theorem~\ref{WeakBound} holds. Indeed, the sign of the matrices do not play any role in the proof, so the probability of $\overline{Y_1}$ can be bounded by the same quantity as in the previous section. In addition, the probability of having a closed walk conditioned to $Y_1$ is also $\frac{1}{n}$ ($n$ is the size of a partition). Therefore, Proposition \ref{WeakBound} applies to bipartite regular graphs.

We can define the analogous recursive relation used in the preceding part. Since this formula represents a quantity that does not depend on the sign of the letters (they are determined by the parity of the positions), we can just ignore them:

$$
Y_p(c,\ell, d) = \left\{
    \begin{array}{ll}
        1 \mbox{ if } c = 0 \mbox{ and } \ell = 0 \\
        0 \mbox{ if } p > \ell \\
        \sum_{i = c}^{\lfloor\frac{\ell}{p}\rfloor} \binom{d}{i}
(p!)^{-i}\frac{\ell!}{(\ell-ip)!}Y_{p+1}(0,\ell-ip, d-i) \mbox{ otherwise.}
    \end{array}
\right.
$$
As before, $Y_p(c,\ell, d)$ is the number of words of size $\ell$ on alphabet of size $d$ that have at least $c$ different letters that appear $p$ times and whose other present letters have more occurrences. The only difference with $X_p(c,\ell, d)$ is that we do not deal with signs. For the same reason, we have

$$
|Y_1| = Y_1(1,2m, d),
$$
$$
|Y_2| = Y_2(1,2m, d),
$$
$$
|Y_3| = Y_3(1,2m, d)
$$
and
$$
|Y_4| = Y_4(0,2m, d).
$$

We have already shown that when one letter appears exactly once, the probability of having a closed walk is $\frac{1}{n}$. Similarly, when a letter appears exactly three times, the probability of getting a closed walk is less than $\frac{5}{n}$. Indeed, in the proof of Lemma \ref{lemmaX3}, all possible configurations of signs for the letter that appears three times are considered (e.g. $\omega = \dots s  \dots s^{-1} \dots s \dots$ or $\dots s \dots s \dots s^{-1} \dots$). No assumption is done on their respective probabilities to occur. These probabilities may or may not be different in the bipartite setting. Since this bound ($\frac{5}{n}$) is the maximum over all the probabilities of getting a closed walk with each configuration of signs, the resulting bound for the probability of getting a closed walk in the bipartite case remains the same.
Hence we get
$$
\mathbf{P}(\omega.v_1 = v_1 | Y_1)\mathbf{P}(Y_1) \le \bigg(\frac{1}{d}\bigg)^{2m}Y_1(1, 2m, d)\frac{1}{n}
$$
and
$$
\mathbf{P}(\omega.v_1 = v_1 | Y_3)\mathbf{P}(Y_3) \le \bigg(\frac{1}{d}\bigg)^{2m}Y_3(1, 2m, d)\frac{5}{n}.
$$
As before, we do not bound the probability of getting a closed walk under condition $Y_4$. Then
$$
\mathbf{P}(\omega.v_1 = v_1 | Y_4)\mathbf{P}(Y_4) \le \bigg(\frac{1}{d}\bigg)^{2m}Y_4(0, 2m, d).
$$

We now estimate $\mathbf{P}(\omega.v_1 = v_1 | Y_2)\mathbf{P}(Y_2)$. We set $y_3(i) = Y_3(1,2m-2i, d-i)$ and $y_4(i) = Y_4(0,2m-2i, d-i)$.

\begin{lemma}

\begin{multline}
\mathbf{P}(\omega.v_1 = v_1 | Y_2)\mathbf{P}(Y_2) \le \\
Y_2(1, 2m, d)\frac{2}{n} + \sum_{i = 1}^m \binom{d}{i} \bigg(\frac{1}{2}\bigg)^i \frac{(2m!)}{(2m-2i)!} \frac{2^i}{(i+1)!}\bigg(y_3(i)\frac{5}{n} + y_4(i) \bigg)
\end{multline}

\end{lemma}

\begin{proof}

We proceed in a similar way as in the proof of lemma \ref{lemmaX2'}. Consider that we have $i$ pairs of matrices that appear exactly  twice in $\omega$. In the bipartite setting, the signs are forced by the parity of the position of each letter. We thus choose to ignore them. Then, the probability of having no pair $(s,t)$ such that we have the pattern $\omega = \dots s \dots t \dots s \dots t \dots$ is $\frac{2^i}{(i+1)!}$. The proof is the same as that of claim 2 in Lemma \ref{lemmaX2'}, except that the numerator and the denominator of the fraction are both divided by $2^i$ (because we ignore the signs).

Using the proof of Lemma \ref{lemmaX2}, we conclude that, if a letter appears twice with the same sign, the probability of having a closed walk is less than $\frac{2}{n}$. If there is a pair $(s, t), \omega = \dots s \dots t \dots s^{-1} \dots t^{-1} \dots$, then, by using the argument from lemma \ref{lemmaX2'} (claim 1), the probability of getting a closed walk is $\frac{1}{n}$.
If those cases do not happen, we still can bound the probability of getting a closed walk using lemma \ref{lemmaX3} when at least a letter appears three times. The conditional probability of getting a closed walk is then less than $\frac{5}{n}$.

Let us combine together all these bounds. 
\begin{itemize}
    \item The union of the event in which a letter appears exactly twice with same sign, and the event where the letters that appear twice form a bad parenthesized word (if we forget about the signs) has size smaller than $Y_2(1, 2m, d)$. The probability of getting a closed walk in this case is not greater than $\frac{2}{n}$.

    \item The size of the event in which this bound does not apply, but we can apply the bound from lemma \ref{lemmaX3} (which is $\frac{5}{n}$) can be computed as follows. $(2i)!\binom{2m}{2i}$ is the number of ways of placing $i$ pairs of letters with different sign in $2m$ positions. Since we ignore the sign, this quantity has to be divided by $2^i$ which gives $(\frac{1}{2})^i \frac{(2m!)}{(2m-2i)!}$. A fraction $\frac{2^i}{(i+1)!}$ of them are well formed. $y_3(i)$ is the number of ways of filling the remaining gaps so that no letter appears once nor twice, and at least letter appears three times. Choosing the $i$ pairs among the $d$ possible ones and summing over all $i \le 2m$, we get
$$
\sum_{i = 1}^m \binom{d}{i} \bigg(\frac{1}{2}\bigg)^i \frac{(2m!)}{(2m-2i)!} \frac{2^i}{(i+1)!}y_3(i).
$$

    \item Similarly, the number of remaining words that correspond to walks whose probability of being closed is not estimated is
$$
\sum_{i = 1}^m \binom{d}{i} \bigg(\frac{1}{2}\bigg)^i \frac{(2m!)}{(2m-2i)!} \frac{2^i}{(i+1)!}y_4(i).
$$

\end{itemize}
Summing all these quantities, multiplying them by their respective probabilities of getting a closed walk and dividing the whole expression by the number of possible words (that is $d^{2m}$), we can conclude.
  \qed \end{proof}


Summing all the above probabilities  and substituting this in \ref{eq:2} finishes the proof of Theorem \ref{BoundBP}.

\subsection{Proof of Proposition \ref{THbireg}}

We now turn to bound the second largest eigenvalue for biregular graphs. Let $G'$ be a bipartite regular graph of degree $d_1$  and whose number of vertices is $2n_1$. Let $\gamma$ be an integer that divides $n_1$. 
We construct the $d_1d_2$-biregular bipartite graph $G$ by merging every $\gamma$ vertices in the right partition. We denote $n_2 = \frac{n_1}{\gamma}$ the size of this partition, and $d_2 = \gamma d_1$, thus $d_1$ and $d_2$ are the respective degrees of each partition of the graph. 
Let 
\[
P = \left(\begin{array}{@{}c|c@{}}
  0 & M \\ \hline \,^t M & 0
\end{array}\right)
\]
be its adjacency matrix. Hence $M$ has dimension $n_1 \times n_2$. Let
\[
Q = \left(\begin{array}{@{}c|c@{}}
  0 & A \\ \hline \,^t A & 0
\end{array}\right)
\]
be the adjacency matrix of $G'$, which is the bipartite regular graph before merging the vertices of the right partition.
We set $J$ such that $M.^tM = A.J.^tA$, thus $J =  I_{n_2} \otimes~J_\gamma$, where $J_\gamma$ is the $\gamma \times \gamma$ matrix whose entries are only ones and $\otimes$ is the Kronecker product. An example of resulting paths is shown in Fig.~\ref{fig:4}.

\begin{figure}[ht]
\begin{center}
\begin{tikzpicture}[roundnode/.style={circle, draw,  minimum width=2pt, inner sep=1pt}]
        \node[roundnode](1){};
        \node[roundnode](2)[below of=1]{};
        \node[roundnode](3)[below of=2]{};
        \node[roundnode](4)[below of=3]{};
        \node[roundnode](5)[below of=4]{};
        \node[roundnode](6)[below of=5]{};
        
        \node[roundnode](a)[right of=1, xshift=1cm]{};
        \node[roundnode](b)[right of=2, xshift=1cm]{};
        \node[roundnode](c)[right of=3, xshift=1cm]{};
        \node[roundnode](d)[right of=4, xshift=1cm]{};
        \node[roundnode](e)[right of=5, xshift=1cm]{};
        \node[roundnode](f)[right of=6, xshift=1cm]{};
        
        \node[roundnode](w)[right of=a]{};
        \node[roundnode](x)[right of=b]{};
        \node[roundnode](y)[right of=c]{};
        \node[roundnode](z)[right of=d]{};
        \node[roundnode](s)[right of=e]{};
        \node[roundnode](t)[right of=f]{};
        
        \node[roundnode](1a)[right of=w, xshift=1cm]{};
        \node[roundnode](1b)[right of=x, xshift=1cm]{};
        \node[roundnode](1c)[right of=y, xshift=1cm]{};
        \node[roundnode](1d)[right of=z, xshift=1cm]{};
        \node[roundnode](1e)[right of=s, xshift=1cm]{};
        \node[roundnode](1f)[right of=t, xshift=1cm]{};
        
        node at (1, 2) {\vdots};

        \draw [->](1) -- (b) node{};
        \draw [->](4) -- (a) node{};
        \draw [->](6) -- (c) node{};
        \draw [->](5) -- (e) node{};
        \draw [->](2) -- (d) node{};
        \draw [->](3) -- (f) node{};
        
        \draw [->](a) -- (w) node{};
        \draw [->](b) -- (w) node{};
        \draw [->](c) -- (w) node{};
        \draw [->](a) -- (x) node{};
        \draw [->](b) -- (x) node{};
        \draw [->](c) -- (x) node{};
        \draw [->](a) -- (y) node{};
        \draw [->](b) -- (y) node{};
        \draw [->](c) -- (y) node{};
        
        \draw [->](d) -- (z) node{};
        \draw [->](e) -- (z) node{};
        \draw [->](f) -- (z) node{};
        \draw [->](d) -- (s) node{};
        \draw [->](e) -- (s) node{};
        \draw [->](f) -- (s) node{};
        \draw [->](d) -- (t) node{};
        \draw [->](e) -- (t) node{};
        \draw [->](f) -- (t) node{};
        
        \draw [->](x) -- (1a) node{};
        \draw [->](w) -- (1d) node{};
        \draw [->](s) -- (1e) node{};
        \draw [->](z) -- (1b) node{};
        \draw [->](y) -- (1f) node{};
        \draw [->](t) -- (1c) node{};
        
        \node(m1)[above of=1, xshift=1cm, yshift=-0.75cm]{$^tA$};
        \node(m2)[above of=a, xshift=0.5cm, yshift=-0.75cm]{$J$};
        \node(m3)[above of=w, xshift=1cm, yshift=-0.75cm]{$A$};

\end{tikzpicture}
\caption{Traveling three steps starting from the left in this directed graph is equivalent to doing two steps in the bipartite graph (starting from the left as well) with merged vertices on the right. The merging operation is represented by $J$ which corresponds to the complete bipartite graphs in the middle. Here $n_1 = 6, n_2 = 2, d_1 = 1, d_2 = 3$.}
\label{fig:4}
\end{center}
\end{figure}

We set $A = (a_{ij})_{i,j \in [\![ 1,\ n_1 ]\!]}$. All $A$'s columns and rows sum up to $d_1$ ---so do those of $^tA$. We can show that for every $x = (x_1, \dots, x_{n_1})$ orthogonal to $e_1 = \frac{1}{\sqrt{n_1}}(1, \dots, 1)$, $Ax$ is also orthogonal to $e$. Indeed, the coordinates of such an $x$ sum up to zero. We denote $Ax = (y_1, \dots, y_{n_1})$. Then

$$
\sum_{i = 1}^{n_1} y_i = \sum_{i = 1}^{n_1} \sum_{j = 1}^{n_1} a_{ij}x_j = \sum_{j = 1}^{n_1} x_j \sum_{i = 1}^{n_1} a_{ij} = d_1 \sum_{j = 1}^{n_1} x_j = 0.
$$
Therefore, $Ax$ is orthogonal to $e$. The same is true for $^tAx$.

On the other hand, it is easy to see that the spectrum of $J$ is
$$
\underbrace{(\gamma, ..., \gamma,}_{n_2} \underbrace{0, ..., 0)}_{n_1 - n_2}.
$$
Let $\lambda_2(X)$ be the second largest eigenvalue of some square matrix $X$ and let $x$ be the normalized eigenvector associated to $\lambda_2(M.^tM)$. For every positive number $q$, we have

$$
|\lambda_2(M.^tM)^q| = ||(A.J.^tA)^q .x|| = ||A(J.^tA.A)^{q-1}J.^tA.x|| \le \gamma d_1 ||A(J.^tA.A)^{q-1} x'||
$$
with $x'$ a normalized vector orthogonal to $e$ (because of the preceding fact). Hence 
$$
||(J.^tA.A)^{q-1} x'|| \le \gamma^{q-1} |\lambda_2^{q-1}(^tA.A)|.
$$
We conclude that 
$$
|\lambda_2(M.^tM)^q| \le d_1^2\gamma^q|\lambda_2(^tA.A)^{q-1}| = d_1d_2(\gamma |\lambda_2(^tA.A)|)^{q-1}.
$$

Since this is a positive quantity, its $q$-th root is defined:
$$
\lambda_2(M.^tM) \le \bigg(d_1d_2(\gamma |\lambda_2(^tA.A)|)^{q-1} \bigg)^{\frac{1}{q}} = \bigg(d_1d_2 \frac{(\gamma|\lambda_2(^tA.A)|)^q}{|\lambda_2(^tA.A)|} \bigg)^{\frac{1}{q}}
$$
which gives
$$
\lambda_2(M.^tM) \le \bigg(\frac{d_1d_2}{|\lambda_2(^tA.A)|}\bigg)^{\frac{1}{q}} \gamma|\lambda_2(^tA.A)|
$$
By taking the limits when $q$ goes to infinity, we obtain
$$
\lambda_2(M.^tM) \le \gamma|\lambda_2(^tA.A)|
$$
We denote by $\mu_1 \ge \mu_2... \ge \mu_{n_1}$ the eigenvalues of $Q$.
To finish the proof

It is well known that for a bipartite graph, the matrix of edges from one partition to the other has singular values that are the squares of that of the adjacency matrix. This way we get that $|\lambda_2(^tA.A)| = \mu_3^2$.
We denote by $\alpha$ some upper bound for $|\mu_3|$. In $^tA.A$, the second largest magnitude  eigenvalue is thus $|\lambda_2(^tA.A)| = (d_1|\mu_3(Q)|)^2 \le d_1^2 \alpha^2$. Hence we have 
$$
|\lambda_2(P)| = \sqrt{|\lambda_2(M.^tM)|} \le \sqrt{\gamma|\lambda_2(A.^t A)|} \le \sqrt{d_1d_2} \alpha.
$$

\subsection{Lower bound for Schreier graphs of Abelian groups}

It is also interesting to test Schreier graphs from Abelian groups. A folklore fact claims that Cayley graphs from such groups cannot be good spectral expanders. This has been shown in \cite{alon-roichman} and later in \cite{abelianSpectralEstimation}. It is interesting to provide an alternative proof that applies also to Schreier graphs. First, as before, we present the experimental results obtained on Cayley graphs from the multiplicative group $\mathbb{Z}/q\mathbb{Z}$ which are actually our Schreier graphs $Sch(GL_k(\mathbb{F}_q) \circlearrowleft (\mathbb{F}_q^k)^*, S)$ with $k = 1$.

\begin{figure}[h!]
\centering
    \includegraphics[trim={50 23 45 41},clip,scale = 0.7]{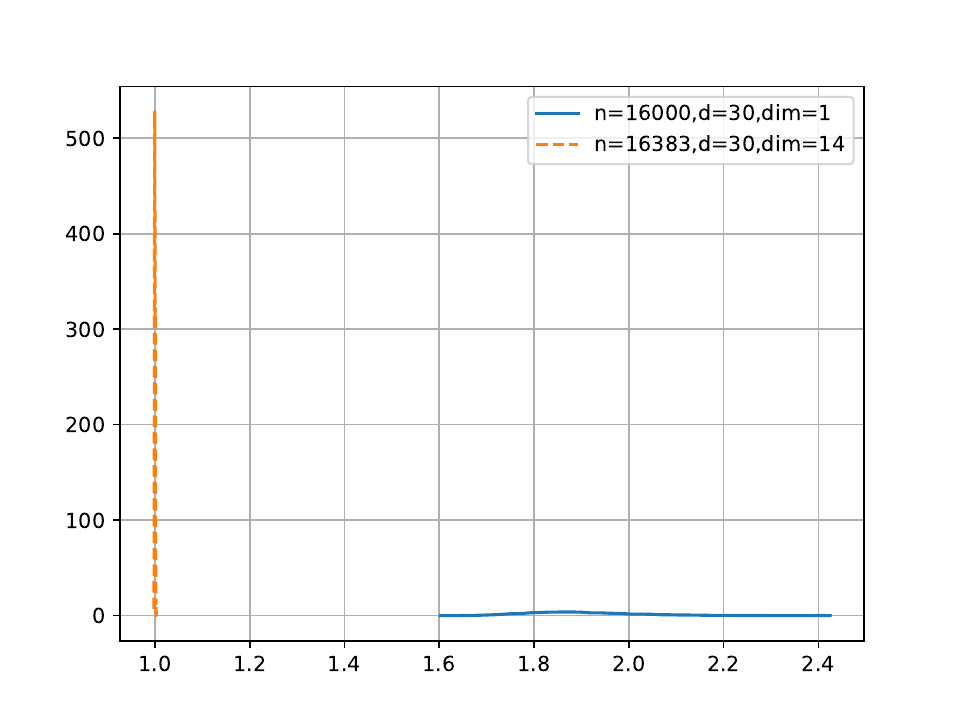}
  \caption{On the left, we display the Schreier graphs from $GL_{14}(\mathbb{F}_2)$ for comparison (the distribution looks very much like that of the permutation model). On the right, we can see how the distribution for Abelian groups is spread out in contrast, and the mean is much larger.}\label{fig:ab-sch}
\end{figure}

  \begin{figure}[!ht]
  \centering
    \includegraphics[trim={50 23 45 41},clip,scale = 0.7]{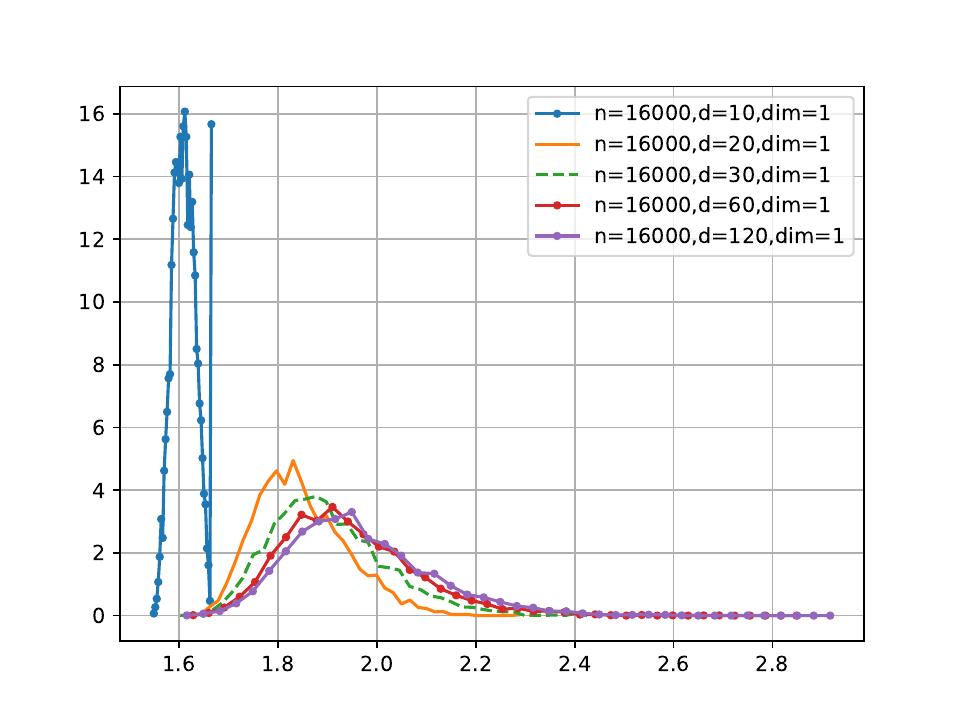}
  \caption{We show the second largest eigenvalue distribution for Cayley graphs from an Abelian group of size $16000$ with different degrees.}\label{fig:abelian}
\end{figure}

In Figure \ref{fig:ab-sch}, we can see the massive difference between nearly optimal spectral expanders and a fairly similar construction with Abelian groups. In comparison, the former second largest eigenvalue distribution seems extremely concentrated around $1$, which, as before, corresponds to the Ramanujan threshold, $2\sqrt{d-1}$. On the other hand, the blue (plain) curve has massive variance and its expected value lies around $1.9 \times 2\sqrt{d-1}$, which makes the associated graphs much worse expanders.

In Figure \ref{fig:abelian}, we can zoom in the distribution of the second largest eigenvalue of such graphs. There is some strange behavior with lower degree graphs, as a significant part of them suddenly get an exceptionally greater second largest eigenvalue than the rest of them. Moreover, it looks significantly smaller in proportion than the other distributions. This is because, since the degree is pretty small, the Ramanujan threshold is quite close to the maximal eigenvalue. The bin on the right of the blue (pointed) curve of $10$-regular graphs must then be the proportion of disconnected graphs whose second largest eigenvalue is $10$. Such phenomenon does not seem to occur with larger degrees, and we observe similar curves as for Schreier graphs of non-Abelian groups, but less concentrated and shifted to the right.

We now can formulate a statement that is an attempt to explain this phenomenon. The crucial point is that in a Schreier graph over an Abelian group, the number of vertices achievable from a fixed vertex $v$ by a path of length $m$, grows polynomially in $m$, while in a good spectral expander this number grows exponentially. We make this observation more precise in the following lemma.

\begin{lemma}\label{l:ABgraphs}
Let $G = (V,E)$, $|V|= n$ be a $2d$-regular Schreier graph from an Abelian group $H$. Let $0 <\eta < 1$. If there exists an integer $m$ such that we have
\[
m \ge -\frac{\log 2n}{2\log\eta}
\]
and
\[
\frac{n}{2} \ge 2^d \binom{m+d}{d},
\]
then, the second largest eigenvalue (in absolute value) of $G$, $|\lambda_2|$, is larger than $\eta 2d$.
\end{lemma}

\begin{proof}
The strategy of the proof is to get a lower bound on the second largest eigenvalue by applying the expander mixing lemma on sets that have no edges between them. Let $v$ be a vertex of $G$ and $A = \{v\}$. Let $B$ be the set of unreached vertices in every random walk of $m$ steps starting at $v$. We will choose $B$ that contains only vertices that cannot be reached with any path of length $m$ starting at $v$ (but not necessarily all such vertices). In order to get a lower bound on the second largest eigenvalue, we want the set $B$ to be non-empty. The exact size of $B$ does not really matter for what comes next (it will only slightly change the choice of the constants), so let us suppose that $|B| = \frac{n}{2}$ for simplicity.

We apply the expander mixing lemma (see \cite{EML} for probably the first proof) to the graph $G^m$ corresponding to the paths of size $m$ with the sets $A$ and $B$.
\[
\bigg||E(A,B)| - (2d)^m\frac{|A||B|}{n} \bigg| \le |\lambda_2|^m \sqrt{|A||B|}.
\]
Since we know the size of $A$ and $B$, and since there are no edges between those sets, ($|E(A,B)| = 0$) the inequality rewrites
\[
(2d)^m \le |\lambda_2|^m \sqrt{2n}.
\]
Let $\mu_2 = \frac{\lambda_2}{2d}$ be the normalized second largest eigenvalue. We get
\[
\bigg(\frac{1}{2n}\bigg)^{\frac{1}{2m}} \le |\mu_2|.
\]
Suppose $0 < |\mu_2| < 1$ (the only graph that has $\mu_2 = 0$ is ``more'' than complete, and if $\mu_2 = 1$, there is nothing to prove). For $|\mu_2|$ to be larger than some $\eta \in ]0,1[$ it is enough to obtain
\[
\eta \le \bigg(\frac{1}{2n}\bigg)^{\frac{1}{2m}} \iff m \ge -\frac{\log 2n}{2\log \eta}.
\]
Which gives the first clause of the assumption in the lemma statement.

It remains to exhibit a condition that is enough  to get $A$ and $B$ such that any walk of $m$ steps from $A$ cannot reach a vertex of $B$. This happens when the number of reachable vertices is smaller than $\frac{n}{2}$. The graph $G$ is defined by the action of $d$ elements of $H$ on each element of $V$. We denote $h_1, \dots, h_d$ those elements, which are chosen at random.

Starting from a vertex $v$, we bound the number of reachable vertices in $m$ steps (denoted by $N_m$) by the number of walks one can do in the graph. We can label each step of the random walk by the element in the group $H$ the traveled edge corresponds to, namely $h_1, \dots, h_d$. For some $i \le d$, $h_i^{-1}$ corresponds to the edge traveled backwards.

Doing so, we have at most $(2d)^m$ random walks. Since $H$ is Abelian, the order of the edges does not matter. Hence we only need to know how many times each $h_i$ appears. Hence, we have at most $\binom{m+2d-1}{m}$ reachable vertices (this is the number of ways of determining how many times each of the $2d$ letters appears in a word of size $m$). This gives the second clause of the assumption of the lemma.
\qed \end{proof}

Lemma \ref{l:ABgraphs} implies in particular that Schreier graphs over an Abelian group with logarithmic degree can be bad spectral expanders. The remaining question is whether or not graphs that verify the conditions of this lemma can be obtained. We claim the following statement:

\begin{prop}\label{ABstatement}
    Let $G$ be a Schreier $2d$-regular graph from an Abelian group on $n$ vertices and let $\lambda_2$ be the second largest eigenvalue of $G$. For every $\epsilon > 0$, there exists $\delta > 0$ such that, when $n$ goes to infinity, 
    \[
    d = \delta \log n \implies |\lambda_2| \ge (1-\epsilon) 2d.
    \]
\end{prop}

\begin{proof}
From the proof of Lemma \ref{l:ABgraphs}, for our set $B$ to be non-empty (in particular, to have size larger than $\frac{n}{2}$), we need the number of reachable vertices to be less than $\frac{n}{2}$. Let us fix $\epsilon$. The bound is trivial if $\epsilon \ge 1$ so let us assume this is not the case. We set $\eta = 1 - \epsilon$. By taking $m$ above our threshold from the preceding lemma
\[
-\frac{\log 2n}{2\log \eta}= \Theta(\log n)
\]
but still logarithmic in $n$ (say $m = \alpha \log n = \Theta(\log n)$ for some constant $\alpha$), we get, by setting $d = O(\log n)$,
\[
N_m \le \binom{m+2d-1}{2d-1} =  \binom{m+2d-1}{m} \le \bigg(\frac{e(m+2d-1)}{m}\bigg)^{2d-1}
\]
By replacing $m$ and $d$ by their order of magnitudes in terms of $n$, 
\[
N_m \le \bigg(\frac{\Theta(\log n)}{\Theta(\log n)}\bigg)^{O(\log n)} = \Theta(1)^{O(\log n)}
\]

Hence, for large $n$, by choosing the constants in the $O$ notation for $d$ (denoted $\delta$ in the statement) and $\alpha$, we can obtain 

\[
N_m \le \frac{n}2.
\]
We thus verify both conditions of Lemma \ref{l:ABgraphs}. Hence, with appropriate choice of the constant $\alpha$ and $\delta$, we can get the desired result.
  \qed \end{proof}

This statement claims that large enough Schreier graphs from an Abelian group with logarithmic degree (with small enough constant in front of the logarithm) can have arbitrarily small spectral gap.

\section{Final comments}

In this paper we study a pseudo-random construction of spectral expanders represented as Schreier graphs. The experimental results suggest that these graphs have nearly optimal value for the second largest eigenvalues (close to $2\sqrt{d-1}$ for regular graphs of degree $d$ and, respectively, to $\sqrt{d_1-1}+\sqrt{d_2-1}$ for biregular bipartite graphs with degrees $d_1$ and $d_2$),
and this is true not asymptotically but for graphs of reasonably small sizes relevant for practical applications. 

A possible direction of the future research is a reduction of the number of random bits used to produce each graph. It would be interesting to reduce the number of these bits from $O(d\log^2 n)$ to $O(d\log n)$, which would bridge the gap between (pseudo)random and deterministic constructions and improve computational complexity of the generation. Our experiments show that a promising direction is a construction based on  Toeplitz matrices.

The ``computer assisted theoretical bounds'' that we provide  are a small step to an explanation of the obtained experimental results. Instead of more traditional asymptotic bounds, we focused on theoretical bounds that can be calculated (possibly with the help of computer) for graphs with specific values of parameters, including the values of parameters that appear  in our numerical experiments.

We would also like to draw attention to the observation that useful theoretical bounds and estimates can be obtained with the help of ``hideous''  formulas combined with computer computations. This approach contradicts the  spirit of classical mathematics, where we appreciate proofs that are elegant and meaningful for a human. However, these techniques  can be useful to justify properties of objects (e.g., of pseudo-random graphs) with specific parameters, which may be necessary in practical applications. Additionally, we provide a simple explanation of the fact that Schreier graphs from Abelian groups can be poor spectral expanders and more data to support this phenomenon.

\section*{Acknowledgments}
This work is part of the author's Ph.D thesis. Special thanks go to his advisers, Andrei Romashchenko and Alexander Shen, for providing several scientific ideas and help in the writing process.


\bibliographystyle{spmpsci}      

\bibliography{main.bib}

\section*{Statements and declarations}

The author was supported in part by the  project ANR FLITTLA ANR-21-CE48-0023.

\noindent The author has no relevant financial or non-financial interests to disclose.

\end{document}